\documentclass[12pt]{article}
\usepackage{graphicx} 

\usepackage{authblk} 
\usepackage{geometry}
\usepackage{amsmath}
\usepackage{amssymb}
\usepackage{esint}
\usepackage{bm}
\usepackage{array}
\usepackage{float}
\usepackage{flafter}
\usepackage{appendix}
\usepackage{caption}
\usepackage{subcaption}
\usepackage{enumitem}
\usepackage{multirow}
\usepackage{float}
\usepackage{rotating}

\usepackage{tikz}
\usetikzlibrary{angles, quotes}

\usepackage[hidelinks]{hyperref}
\hypersetup{colorlinks=true,
            citecolor=blue,}

\newcommand{\pfrac}[2]{\frac{\partial{#1}}{\partial{#2}}}

\newcommand{\lbrac}[0]{\left(}
\newcommand{\rbrac}[0]{\right)}
\newcommand{\bracked}[1]{\left( #1 \right)}
\newcommand{\bracfrac}[2]{\left( \frac{#1}{#2} \right)}

\usepackage{xcolor}

\usepackage[style = apa, uniquename=false]{biblatex}
\DeclareLanguageMapping{english}{english-apa}
\title{Stability analyses of divergence and vorticity damping on gnomonic cubed-sphere grids}
\author[1]{Timothy C. Andrews \thanks{timand@umich.edu}}
\author[1]{Christiane Jablonowski}
\affil[1]{Department of Climate and Space Sciences and Engineering, University of Michigan, Ann Arbor, MI, USA}

\begin{document}

\maketitle

\begin{abstract}
Divergence and vorticity damping, which operate upon horizontal divergence and relative vorticity, are explicit diffusion mechanisms used in dynamical cores to ensure stability. To avoid numerical blow-up from excessively strong diffusion, there are mesh-dependent upper bounds on the coefficients of the diffusion operators. This work considers such stability limits for three gnomonic cubed-sphere meshes: the 1) equidistant, 2) equiangular, and 3) equi-edge mappings. Stability limits are derived from a von Neumann analysis of damping with a simplified pseudo-Laplacian operator, as used in NOAA GFDL's finite-volume dynamical core on the cubed-sphere (FV3), and with the full curvilinear Laplacian. The resulting stability limits depend on the gnomonic mapping through the cubed-sphere cell areas, aspect ratios, and grid nonorthogonality. The analytical stability limits are compared to practical divergence and vorticity damping upper bounds in FV3, using idealised tests and the equiangular and equi-edge grids. For divergence damping, both the magnitude of maximum stable coefficients and the locations of instability agree with linear theory. Due to implicit vorticity diffusion in the FV3 transport scheme, practical limits for vorticity damping are lower than the explicit stability limits and depend on the choice of horizontal transport scheme.
\end{abstract}

\section{Introduction}
\label{sec:intro}
Numerical diffusion is an important feature of dynamical cores for weather and climate modelling. It can be implemented in various ways to suit the target application and model design \parencite{jablonowski2011pros_diffusion_chapter}, with all techniques sharing the primary objective of damping small-scale phenomena that could otherwise cause numerical instability. A major reason for this instability is the discrete nature of numerical models; they can only resolve a finite number of wave frequencies, the most oscillatory of these with a wavelength of $2\Delta x$, where $\Delta x$ denotes the grid spacing. Waves near the $2\Delta x$ scale can accumulate an unphysical amount of energy if the model does not properly transfer it to smaller, unresolved scales. The inability of energy to transition to smaller waves is often termed spectral blocking \parencite{boyd2001chebyshev}, and if left unchecked, may result in a numerical simulation blowing up near the grid-scale. To avoid instability, all dynamical cores contain some form of dissipation, such as of the energy or the potential enstrophy. \par
Numerical diffusion can be introduced implicitly through the numerical methods, such as the transport or time stepping schemes, or applied explicitly through additional terms in the dynamical equations. An example of explicit diffusion is the addition of a Laplacian or higher-order operator acting on the velocity components in the horizontal momentum equation. An alternative that this paper considers is divergence and vorticity damping, which instead act upon the horizontal divergence and relative vorticity generated by the horizontal velocity field.  \par
The impact of divergence and vorticity damping in a dynamical core depends on two key parameters: the damping coefficient, which governs the diffusion strength, and the order, which dictates the range of wavenumbers that are appreciably damped. For example, \textcite{carley2023mitigation} found that in the NOAA Geophysical Fluid Dynamics Laboratory (GFDL) finite volume cubed-sphere dynamical core of FV3 \parencite{harris2021scientific_FV3_doc,lin1997explicit,lin2004vertically}, moving from sixth- to eighth-order divergence damping improved performance in a moist squall line test. Conversely, our FV3 studies with a baroclinic wave test encounter instability with eighth-order divergence damping and certain transport schemes, but this instability is avoided when switching to sixth-order damping. Divergence damping parameters have also been noted to affect the modelling of physical phenomena, such as radiative-convective equilibrium \parencite{anber2018sensitivity} and the generation of tropical cyclones and hurricanes \parencite{zhao2012some}. \par
An important consideration for selecting the divergence and vorticity damping coefficients is the existence of both lower and upper bounds for stability with explicit time stepping methods. Insufficient diffusion allows too much kinetic energy and enstrophy at small scales and leads to numerical blow-up. On the other hand, excessive diffusion will exceed a linear stability restriction, also causing a blow-up. The upper bound depends on the computational grid, and this work will investigate the impact of different gnomonic cubed-sphere grids on the damping coefficient limits. \par
With the advent of new computer architectures, many recent and in-development weather and climate models are selecting grids that can exploit the potential for greater spatial parallelism, thus reducing wall-clock times. The gnomonic cubed-sphere is one such competitive grid, and is used in the GungHo dynamical core from the UK Met Office's LFRic-Atmosphere model \parencite{adams2019lfric,melvin2019mixed,melvin2024mixed}, the spectral element (SE) dynamical core \parencite{dennis2012CAM_SE,lauritzen2018ncar_CAM_SE,taylor2010compatible,taylor1997spectral} in NCAR's Community Atmosphere Model (CAM) \parencite{neale2010description}, and FV3. LFRic and CAM-SE use the most common equiangular gnomonic cubed-sphere, whilst FV3 employs a unique equi-edge variant for improved cell regularity along the panel edges. These two grids, along with the original equidistant cubed-sphere of \textcite{sadourny1972_equidistant}, will be examined in this paper. \par
Default diffusion coefficients in dynamical cores are generally tuned through experimentation and should only be considered as guidelines. Modifications to the diffusion strength are often necessary, and knowledge of analytical upper bounds allows diffusive instabilities to be avoided with such changes. We will use the von Neumann method to derive stability limits for divergence and vorticity damping. Although a simple technique, von Neumann analyses can lead to accurate and useful estimations of stability bounds; examples in dynamical cores include for transport schemes in \textcite{kent2014determining} and \textcite{lauritzen2007stability}, and for second- and fourth-order divergence damping on longitude-latitude grids in \textcite{whitehead2011stability}. Whilst this work has similarities to \textcite{whitehead2011stability}, our analysis instead considers cubed-sphere grids and an arbitrary order of divergence damping. Additionally, we investigate the stability of vorticity damping, which is less discussed in the literature than divergence damping. This is because the diffusion coefficient is typically larger for divergence than vorticity damping --- divergent modes are often considered noisy and require filtering, whilst rotational modes are a crucial component of large-scale atmospheric flows. \par
Our linear stability analysis focuses on the staggered Arakawa C- and D-grids \parencite{arakawa1977computational}, where linear stability limits for divergence and vorticity damping are nearly equal, but have small differences due to the offset in divergence and vorticity locations. We will first analyse damping operators that use a simplified version of the full curvilinear Laplacian, as is the case in FV3. We then compare this to the stability of using the full Laplacian. Afterwards, the theoretical stability limits are compared to practical damping limits in the hydrostatic FV3 dynamical core in CAM (CAM-FV3), which supports both the equiangular and equi-edge grids. FV3's transport scheme introduces implicit diffusion of the relative vorticity but does not directly affect the horizontal divergence  \parencite{harris2021scientific_FV3_doc}. Explicit divergence damping is thus necessary for stability, especially given that the $2 \Delta x$ inertia-gravity modes are stationary on the FV3 D-grid \parencite{skamarock2008linear}, whilst vorticity damping is often optional. As will be seen from CAM-FV3 simulations, the allowable magnitude of divergence damping can be accurately predicted by linear stability theory, whilst the practical limit on explicit vorticity damping is much lower due to the implicit diffusion. We hypothesise that the practical vorticity damping limits may reflect the level of implicit transport diffusion in a specific flow. \par
We now outline the content of this paper: Section \ref{sec:prelim} provides background material for divergence and vorticity damping and the staggered Arakawa grids. Section \ref{sec:cubed_sphere_grids} overviews the gnomonic cubed-sphere and compares three different mappings (equidistant, equiangular, equi-edge). The von Neumann stability analysis for divergence and vorticity damping is described in Section \ref{section:linear_stab_analysis}, along with comparisons of stable coefficients and damping of different wavenumbers on the three grids. Section \ref{sec:fv3_and_test_case} compares the linear theory to maximum stable coefficients in baroclinic wave and Held-Suarez simulations with CAM-FV3. Lastly, Section \ref{section:discussion} will discuss the key findings and future research directions. \par

\section{Background material }
\label{sec:prelim}

\subsection{Divergence and vorticity damping}
\label{sec:div_vort_damp}
In the horizontal momentum equations, an explicit diffusion of order $2q, q \in \mathbb{Z}^{+},$ can be applied directly to the prognostic variables through
\begin{equation}
    \pfrac{\mathbf{u}}{t} = ... + (-1)^{q+1} \nu \bm{\nabla}^{2q} \mathbf{u},
\label{eq:velocity_damp}
\end{equation}

\noindent where $\bm{\nabla}^{2}$ is the vector Laplacian containing only horizontal derivatives, $\mathbf{u} = (u,v,0)$ symbolises the horizontal velocity field, with $u$ and $v$ the zonal and meridional components of the wind, and $\nu$ denotes the diffusion coefficient. Laplacian diffusion, which represents physical viscosity in the Navier-Stokes equations, is obtained when $q=1$, whilst $q \geq 2$ are hyperviscosities which are less diffusive for larger-scale waves. Fourth-order damping ($q=2$) is the most common choice in dynamical cores \parencite{jablonowski2011pros_diffusion_chapter}. \par
Divergence and vorticity damping provide an alternative to horizontal velocity damping (\ref{eq:velocity_damp}), where energy is instead removed from the divergent and rotational modes. This exploits the decomposition of the Laplacian into
\begin{equation}
    \bm{\nabla}^{2} \mathbf{u} = \nabla D + \nabla \times \zeta \hat{\mathbf{k}},
\end{equation}

\noindent with scalar quantities of the horizontal divergence, $D = \nabla \cdot \mathbf{u}$, and relative vorticity in the horizontal plane, $\zeta = \hat{\mathbf{k}} \cdot (\nabla \times \mathbf{u})$, where $\hat{\mathbf{k}} = (0,0,1)$ is the unit vector in the vertical direction. Introducing these quantities into the horizontal momentum equations results in

\begin{equation}
    \pfrac{\mathbf{u}}{t} = ... + (-1)^{q+1} \nu_D \nabla \bracked{\nabla^{2(q - 1)} D} + (-1)^{q+1} \nu_{\zeta} \nabla \times \bracked{\nabla^{2(q - 1)} \zeta \mathbf{\hat{k}}},
\label{eq:div_vort_damp}
\end{equation}

\noindent with $\nu_D$, $\nu_{\zeta}$ the coefficients of the divergence and vorticity damping. Computing the divergence and the vertical component of the curl of (\ref{eq:div_vort_damp}) shows that
\begin{subequations}
\begin{align}
    \pfrac{D}{t} &= ... + (-1)^{q+1} \nu_D \nabla^{2q} D, \label{eq:div_time_ev} \\
    \pfrac{\zeta}{t} &= ... + (-1)^{q+1} \nu_{\zeta} \nabla^{2q} \zeta. \label{eq:vort_time_ev}
\end{align}
\label{eq:time_ev_eqs}
\end{subequations}

\noindent Hence, (\ref{eq:div_vort_damp}) leads to a Laplacian or higher-order filtering of the divergent and rotational modes. The damping of each mode is independent \parencite{shuman1969currently,mcpherson1973noise}, so the choice of $\nu_D$ or $\nu_{\zeta}$ does not constrain the other, and in general, $\nu_D \neq \nu_{\zeta}$. Note that the divergence damping in (\ref{eq:div_vort_damp}), which damps two-dimensional horizontal motions, should not be confused with three-dimensional divergence damping of acoustic modes \parencite{skamarock1992stability,klemp2018damping,jablonowski2011pros_diffusion_chapter}. \par
Section \ref{section:linear_stab_analysis} will perform a von Neumann stability analysis on the $\nabla^{2q} s$ operator, for $s \in \{D, \zeta \}$. These limits are presented in terms of a nondimensional damping coefficient of $C_{s, 2q}$, which is defined through
\begin{equation}
    \nu_{s, 2q} = \frac{(C_{s, 2q} ~\Delta A_{\text{min}})^q}{\Delta t},
\label{eq:nonD_damp_coeff}
\end{equation}

\noindent with $\Delta A_{\text{min}}$ the minimum cell area over the cubed-sphere and $\Delta t$ the time step size. The benefit of working with a nondimensional coefficient is that it can be set independently of the horizontal grid resolution. Note that using the minimum area in (\ref{eq:nonD_damp_coeff}) means that $\nu_{s, 2q}$ is the same at each grid cell. Another choice is to use the individual cell areas $\Delta A$ instead of $\Delta A_{\text{min}}$, which causes $\nu_{s, 2q}$ to vary between cells. \par

\subsection{Horizontal grids}
This work will focus on the Arakawa C- and D-grids, which stagger the horizontal wind components. On the C-grid, horizontal divergence is evaluated at cell centres and relative vorticity at cell corners (Fig. \ref{fig:C_grid}) whilst these are swapped on the D-grid, with the vorticity at cell centres and divergence at cell corners (Fig. \ref{fig:D_grid}). A D-grid can be constructed on the C-grid by using the cell centres as vertices, so that vorticity is at cell centres and divergence at cell corners; likewise, a C-grid can be formed by connecting D-grid cell centres (Fig. \ref{fig:C_D_grids}). We will refer to this second grid as the \textit{offset} grid, with cell centre quantities on the primary grid being cell corner quantities on the offset grid, and vice versa. \par
For a resolution with $N_e \times N_e$ cell centre values on each panel of the primary grid, there will be $(N_e + 1) \times (N_e + 1)$ values on the offset grid panel. Cell centre values on the primary grid will be denoted by integer indices, $s_{i,j}$, with $ ~i,j \in \{1, 2, ...,N_e \}$ denoting the column and row indices. Indices for cell centre values on the offset grid are offset by one-half, $s_{i - 0.5, j - 0.5}, ~i,j \in \{1, 2, ..., N_e+1 \}$. Linear stability limits for divergence or vorticity damping are applied on the grid (primary or offset) where that quantity is at the cell centres. For example, with the D-grid, the offset grid is used for divergence damping and the primary grid for vorticity damping. Our analysis can also be applied to other horizontal grids, given a structured layout of the horizontal divergence and relative vorticity on each cubed-sphere panel. \par

\begin{figure}[htpb]
    \centering
    \begin{subfigure}{0.49\textwidth}
    \begin{tikzpicture}[scale=2.5]
    \draw[step=1cm, gray, thin, dashed] (0,0) grid (3,3);
    \foreach \x in {0.5,1.5,2.5}
        \foreach \y in {0.5,1.5,2.5}
            \fill (\x,\y) circle[radius=1pt] node[above right] {$\zeta_{\x,\y}$};
    \foreach \x in {1,2}
        \foreach \y in {1,2}
            \fill (\x,\y) circle[radius=1pt] node[above right] {$D_{\x,\y}$};

    \foreach \x in {0.5,1.5,2.5}
        \foreach \y in {0.5, 1.5}
            \draw[solid] (\x, \y) -- (\x, \y+1);

    \foreach \y in {0.5,1.5,2.5}
        \foreach \x in {0.5, 1.5}
            \draw[solid] (\x, \y) -- (\x+1, \y);

    \draw[red, ultra thick, ->] (0.2,2.0)   -- (0.8,2.0) node[below left] {$u$};
    \draw[red, ultra thick,->] (1.2,2.0)   -- (1.8,2.0) node[below left] {$u$};
    \draw[red, ultra thick,->] (1.0,1.2)   -- (1.0,1.8) node[below right] {$v$};
    \draw[red, ultra thick,->] (1.0,2.2)   -- (1.0,2.8) node[below right] {$v$};
        
    \end{tikzpicture}
    
    \caption{C-grid}
    \label{fig:C_grid}
    \end{subfigure}
    \hfill
    \begin{subfigure}{0.49\textwidth}
        \begin{tikzpicture}[scale=2.5]
    \draw[step=1cm, gray, thin, dashed] (0,0) grid (3,3);
    \foreach \x in {0.5,1.5,2.5}
        \foreach \y in {0.5,1.5,2.5}
            \fill (\x,\y) circle[radius=1pt] node[above right] {$D_{\x,\y}$};
    \foreach \x in {1,2}
        \foreach \y in {1,2}
            \fill (\x,\y) circle[radius=1pt] node[above right] {$\zeta_{\x,\y}$};

    \foreach \x in {0.5,1.5,2.5}
        \foreach \y in {0.5, 1.5}
            \draw[solid] (\x, \y) -- (\x, \y+1);

    \foreach \y in {0.5,1.5,2.5}
        \foreach \x in {0.5, 1.5}
            \draw[solid] (\x, \y) -- (\x+1, \y);

    \draw[red, ultra thick, ->] (0.5,1.7)   -- (0.5,2.3) node[below right] {$v$};
    \draw[red, ultra thick, ->] (1.5,1.7)   -- (1.5,2.3) node[below right] {$v$};
    \draw[red, ultra thick, ->] (0.7,1.5)   -- (1.3,1.5) node[below left] {$u$};
    \draw[red, ultra thick, ->] (0.7,2.5)   -- (1.3,2.5) node[below left] {$u$};
    \end{tikzpicture}
    
    \caption{D-grid}
    \label{fig:D_grid}
    \end{subfigure}
    \caption[]{Diagrams of the Arakawa C- and D-grids, showing the horizontal wind ($u,v$) components (red) in the top-left cell. The C-grid stores normal wind components on the cell edges, leading to divergence at the cell centre and vorticity at the cell corners. The D-grid stores tangential wind components on the cell edges, leading to vorticity at the cell centre and divergence at the cell corners. The primary grid is shown in solid black, with the offset grid in dashed lines formed by connecting cell centres on the primary grid; an offset D-grid is generated on the C-grid, and an offset C-grid on the D-grid.}
    \label{fig:C_D_grids}
\end{figure}

\section{Gnomonic cubed-sphere grids}
\label{sec:cubed_sphere_grids}
\subsection{Description}
Many dynamical cores in the previous decades used a longitude-latitude (lon-lat) grid, as it is structured, orthogonal, and intuitive for interpreting weather and climate forecasts \parencite{williamson2007evolution}. However, cell areas drastically reduce near the polar singularities on lon-lat grids, which leads to severe time step restrictions or the need for filtering to ensure stability, such as through polar Fourier filters \parencite{umscheid1971further,williamson1973comparison}. On modern computer architectures, the small polar cells greatly reduce the potential wall-clock time improvements available with GPUs and parallelism. An alternative with better scalability is quasi-uniform grids based on refined polyhedra \parencite{staniforth2012horizontal}. One choice is the cubed-sphere, which projects six square panels, corresponding to the faces of a cube, onto the spherical domain. \par
Instances of the cubed-sphere can be broadly categorised as either conformal or gnomonic \parencite{staniforth2012horizontal}, with possible modifications such as refinement using spring dynamics \parencite{tomita2001shallow}, the use of an elliptic solver \parencite{putman_lin_2007}, or the enforcement of a uniform Jacobian \parencite{ranvcic2017nonhydrostatic}. The conformal mapping \parencite{ranvcic1996global,mcgregor1996semi,adcroft2004implementation} equates angles on the reference cube and the sphere, allowing for orthogonality in the final grid, barring at eight singularities corresponding to the reference cube corners. The gnomonic projection maps straight lines from the cube faces onto great circles on the sphere, which leads to a nonorthogonal grid, but more uniform cell areas compared to a standard conformal mapping. Whilst the coordinates of each panel on the gnomonic cubed-sphere are free of singularities, there are discontinuities at the panel edges. We will focus on the gnomonic variant of the cubed-sphere, which was identified as effective for finite-volume transport in \textcite{putman_lin_2007}.  \par
We will compare three gnomonic cubed-sphere grids, although other variants exist, e.g. \textcite{ purser2017minor} and \textcite{purser2018mobius}. The first is the original \textit{equidistant} grid of \textcite{sadourny1972_equidistant}, which defines cells of uniform length over the reference cube faces. However, this leads to a large range of cell areas when projecting onto the cubed-sphere panels. The second grid is the \textit{equiangular} mapping proposed by \textcite{ronchi1996_equiangular}. This prioritises the uniformity of cells on the cubed-sphere by defining panel coordinates that correspond to uniform angular spacings. The reduced range of equiangular cell areas can greatly improve solution accuracy compared to the equidistant grid, e.g. \textcite{nair2005discontinuous}. The equiangular mapping is used by the GungHo (LFRic) and SE dynamical cores. The third gnomonic grid we consider is the \textit{equi-edge} grid, which was first documented in \textcite{chen2021lmars} and is unique to FV3. The equi-edge mapping is a modification of the equiangular projection with greater uniformity at the panel edges, with the motivation being to reduce grid imprinting errors. \par
Using the gnomonic cubed-sphere requires careful treatment of the panel edge discontinuities, particularly for operations that require information from multiple panels, such as transport, or divergence and vorticity damping at or near the panel edges. In these instances, values from adjacent panels are interpolated onto `ghost' cells that form an extension of the panel of interest; in FV3, the duo grid of \textcite{mouallem2023duo_grid} is used for more accurate ghost cell remapping. Operations that use a larger stencil, such as higher-order hyperviscosities, require more layers of ghost cells; $q$ layers of ghost cells are necessary for damping of order $2q$. Our analyses will assume that ghost cell information is present and will disregard interpolation errors introduced in this process. \par

\subsection{Details and comparisons}
\label{subsec:cube_sphere_grids_descript}
We now discuss the construction of a gnomonic cubed-sphere grid and the difference between the equidistant, equiangular, and equi-edge mappings. Further discussion and visualisations of these mappings can be found in \textcite{chen2021lmars} and \textcite{santos2024analysis_phd}. \par
First, the spatial resolution of the reference cube is chosen and denoted as C$N_e$, with $N_e$ the number of cells along each cubed-sphere panel edge. Next, a range of $N_e + 1$ reference angles, $\theta$, are constructed as
\begin{equation}
    \theta \in (-\theta_{\text{max}}, \theta_{\text{max}}), ~\Delta \theta = \frac{2 \theta_{\text{max}}}{N_e}.
\end{equation}

\noindent Then, local Euclidean coordinates of ($\xi,\eta$) on each cube face are defined as
\begin{equation}
    \xi(\theta) = a \beta(\theta), ~\eta(\theta) = a \beta(\theta),
\label{eq:gnomonic_proj}
\end{equation}
\noindent with $a = R/\sqrt{3}$ and $R$ denoting the Earth's radius. The $\beta(\theta)$ mapping and maximum panel angle $\theta_{\text{max}}$ for the three gnomonic grids are provided in Table \ref{table:cubed_sphere_key_quantities}. We then denote $(x,y)$ as the gnomonic projection of the $(\xi,\eta)$ cube face coordinates onto the surface of the sphere. The curvilinear $(x,y)$ coordinates are free of singularities, but are discontinuous between panels. Additionally, one can describe the gnomonic projection of the reference cube onto three-dimensional Cartesian coordinates using
\begin{equation}
    (X,Y,Z) = \frac{R}{\sqrt{a^2 + \xi^2 + \eta^2}} \mathbf{P}_i,
\label{eq:panel_to_3d_cart}
\end{equation}

\noindent where $\mathbf{P}_i$ denotes a three-dimensional Cartesian representation of the two-dimensional local coordinates for cube panel $i$, $i \in \{ 1,2,3,4,5,6 \}$. For example, the first panel can be defined as $\mathbf{P}_1 = (a,\xi,\eta)$. \par

\renewcommand\arraystretch{2}
\begin{table}[htpb]
\begin{center}
\caption[]{A table of key properties for the gnomonic cubed-sphere grids. The first two rows give parameters for defining the cubed-sphere panel coordinates. The remaining rows show important measurements that are computed numerically to three decimal places at a C192 resolution. Cell areas are computed on the primary grid, so the equiangular ratio of maximum to minimum cell is not quite $\sqrt{2}$. Cell aspect ratios of $\chi = \Delta y/\Delta x$ and the nonorthogonality metric term of $\sin(\alpha)$ are computed on the offset grid, where the cell centre variable lies exactly at the panel corners and the middle of the edges. Note that $\sqrt{2}\approx1.414$ and $\sqrt{3}/2 \approx 0.866$.}
\label{table:cubed_sphere_key_quantities}
\begin{tabular}{ |c||c|c|c| } 
\hline
 & Equidistant & Equiangular & Equi-edge \\
\hline\hline
$\theta_{\max}$ & 1 & $\frac{\pi}{4}$ & $\arcsin\bracfrac{1}{\sqrt{3}}$ \\
\hline
$\beta(\theta)$ & $\theta$ & $\tan(\theta)$ & $\sqrt{2} \tan(\theta)$ \\
\hline\hline
Ratio of max/min cell area & 5.142 & 1.408 & 2.299 \\
\hline
Location of smallest cell & Corners & Middle of edges & Corners \\
\hline
$\chi$ at panel corners & 1.000 & 1.000 & 1.000 \\ 
\hline
$\chi$ at middle of panel edges & 1.414, 0.707 & 1.414, 0.707 & 1.061, 0.943 \\ 
\hline
$\sin(\alpha)$ at panel corners & 0.866 & 0.866 & 0.866 \\
\hline
$\sin(\alpha)$ at middle of panel edges & 1.000 & 1.000 & 1.000 \\
\hline
\end{tabular}
\end{center}
\end{table}

\renewcommand\arraystretch{1}

We compute distances between points on the sphere by using three-dimensional Cartesian coordinates (\ref{eq:panel_to_3d_cart}), then converting these to longitude ($\lambda$) and latitude ($\phi$) coordinates, with
\begin{equation}
    \lambda = \arctan \bracfrac{Y}{X}, \phi=\arcsin \bracfrac{Z}{\sqrt{X^2 + Y^2 + Z^2}}.
\label{eq:cart_to_lonlat}
\end{equation}

\noindent The distance between two longitude-latitude points of $(\lambda_1, \phi_1)$ and $(\lambda_2, \phi_2)$ is then given by the great circle distance of
\begin{equation}
    \Delta d_{12} = R \arccos ( \sin(\phi_1) \sin(\phi_2) + \cos(\phi_1) \cos(\phi_2) \cos(\lambda_1 - \lambda_2)).
\label{eq:great_circle_dist}
\end{equation}

\noindent This procedure is used to compute distances of $\Delta x$ and $\Delta y$ between vertices on the cubed-sphere. \par
To compute cell areas on the cubed-sphere, we require the internal angles, $\alpha$, of each cell. Consider a cell with corner points specified by three-dimensional Cartesian vectors of $\mathbf{p}_1, \mathbf{p}_2, \mathbf{p}_3, \mathbf{p}_4$. As in Appendix B of \textcite{chen2021lmars}, we construct unit vectors of
\begin{equation}
    \hat{\mathbf{e}}_{ab} = \frac{\mathbf{p}_a \times \mathbf{p}_b}{||\mathbf{p}_a \times \mathbf{p}_b||}.
\label{eq:noncoord_basis_vectors}
\end{equation}

\noindent such that at vertex $b$, which has adjacent vertices $a$ and $c$, the interior angle is approximated as
\begin{equation}
    \alpha_{abc} = \arccos(\hat{\mathbf{e}}_{ba} \cdot \hat{\mathbf{e}}_{bc}).
    \label{eq:alpha_comp}
\end{equation}

\noindent The cell area can be computed numerically using the spherical excess formula applied to a quadrilateral \parencite{adcroft2004implementation},
\begin{equation}
    \Delta A = R^2 (\alpha_{412} + \alpha_{123} + \alpha_{234} + \alpha_{341} - 2 \pi).
\label{eq:spherical_excess_area}
\end{equation}

We approximate the cell centre $\alpha$ as the mean of the four corner angles,
\begin{equation}
    \alpha = \frac{1}{4} (\alpha_{412} + \alpha_{123} + \alpha_{234} + \alpha_{341}).
\label{eq:alpha_cell}
\end{equation}

\noindent The cell centre $\alpha$ is used in the $\sin(\alpha)$ metric terms that arise from nonorthogonality of gnomonic cubed-sphere grids and will be present in the later stability analysis. These metric terms are absent in orthogonal coordinate systems where $\alpha = \pi/2$ at all points, so $\sin(\alpha)=1$. The gnomonic cubed-sphere is orthogonal at central lines of $x=0$ and $y=0$ on each panel, with increasing nonorthogonality moving towards the panel corners, where the intersection of three panels leads to $\alpha = 2\pi/3$ and $\sin(\alpha) = \sqrt{3}/2$. \par
\begin{figure}[htpb]
    \centering
     \includegraphics[width=\textwidth]{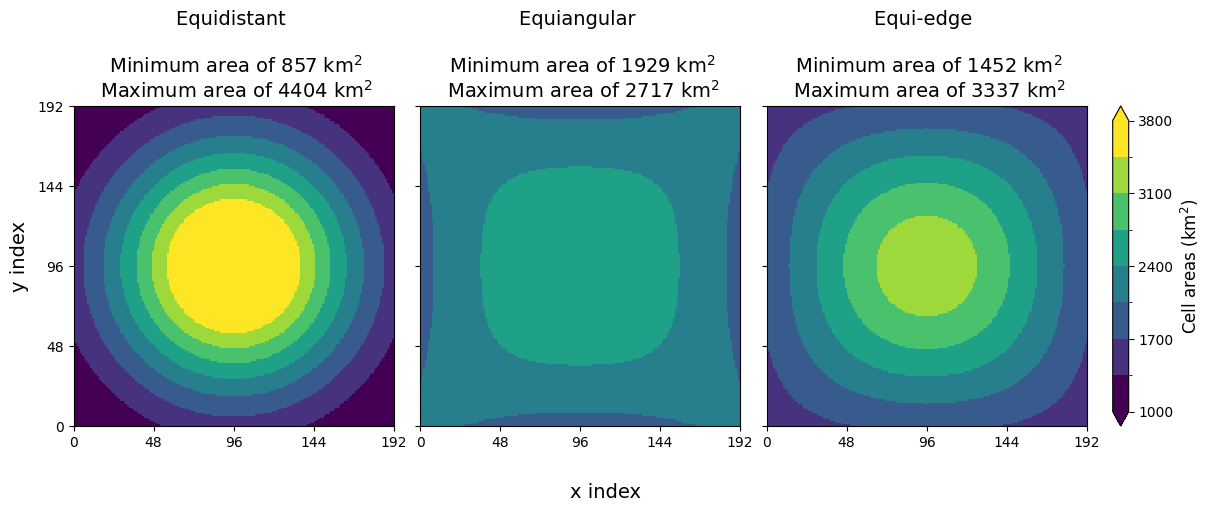}
    \caption[]{Cell areas, computed numerically with (\ref{eq:spherical_excess_area}), for the cubed-sphere (primary) grids at a C192 resolution. The equidistant grid has the largest variation in areas, and the equiangular grid has the narrowest. The smallest cells for the equidistant and equi-edge grids are at the panel corners, whereas the smallest equiangular cells are at the middle of the panel edges.}
    \label{fig:cubed_sphere_areas}
\end{figure}

Table \ref{table:cubed_sphere_key_quantities} compares grid properties that will be important for the diffusive stability limits, specifically the range and distribution of cell areas, and the cell aspect ratios. Cell areas, shown for each mapping at a C192 resolution in Fig. \ref{fig:cubed_sphere_areas}, are largest at the centre of the panels. The equiangular grid has the narrowest range of cell areas, and when $N_e$ is odd, the ratio of the maximum to minimum cell area approaches $\sqrt{2}$ as the grid is refined \parencite{ronchi1996_equiangular}. The equidistant grid has the largest range of cell areas, with its ratio of largest to smallest cell area more than three times that of the equiangular grid (Table \ref{table:cubed_sphere_key_quantities}). The smallest cell for both the equidistant and equi-edge grids is at the panel corners, but is at the middle of panel edges for the equiangular grid. \par
\begin{figure}[htpb]
    \centering
     \includegraphics[width=\textwidth]{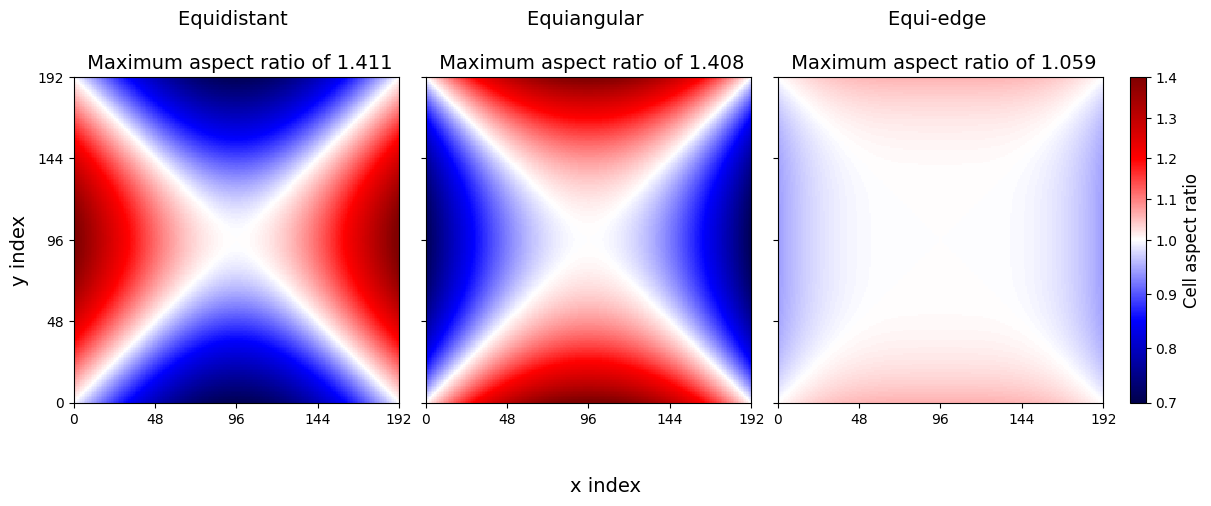}
    \caption[]{Cell aspect ratios of $\chi = \Delta y/ \Delta x$ for the cubed-sphere (primary) grids at a C192 resolution. The largest/smallest aspect ratios are at the middle of the panel edges for all grids. As the equi-edge grid prioritises more uniform cells along the panel edges, it has the narrowest range of aspect ratios. The maximum aspect ratio for the equidistant and equiangular grids is close to $\sqrt{2}$, but not exactly so, due to the even number of cells on the primary grid, i.e. the middle cells are offset from the centre by $\Delta x/2$ or $\Delta y/2$. }
    \label{fig:cubed_sphere_chi}
\end{figure}

Another important property is the cell aspect ratio of $\chi = \Delta y/ \Delta x$. These are shown for the three mappings in Fig. \ref{fig:cubed_sphere_chi}, with all grids containing unity aspect ratios at the panel corners and along the diagonals, with monotonically increasing or decreasing $\chi$ moving from the diagonals towards the central lines of the panel ($x=0$ and $y=0$). The equi-edge grid intentionally has the most uniform aspect ratios of the three grids, with the trade-off being a wider range of cell areas than the equiangular grid. The aspect ratios are similar for the equidistant and equiangular grids, with a maximum value of $\sqrt{2}$ at the middle of the panel edges. \par 

\subsection{Divergence, vorticity, and the Laplacian operator}
\label{subsec:operators}
The definitions of the divergence, vorticity, and Laplacian operators in a general curvilinear coordinate require mapping factors that are encapsulated in the metric tensor. One approach is to form the metric tensor using the coordinate basis, as is computed for the equiangular mapping in \textcite{ronchi1996_equiangular} and \textcite{nair2005discontinuous}. However, the use of the coordinate basis leads to a different metric tensor for each gnomonic mapping and expressions that can be unwieldy to work with. Instead, we use the noncoordinate basis approach of \textcite{chen2021lmars}, which uses the numerically constructed local basis vectors of $\bf{\hat{e}}_1, \bf{\hat{e}}_2$ (\ref{eq:noncoord_basis_vectors}), with angle between them $\alpha = \arccos(\bf{\hat{e}}_1 \cdot \bf{\hat{e}}_2)$. This leads to a simple metric tensor that has the same mathematical form for any gnomonic mapping, allowing for simultaneous analyses of the equidistant, equiangular, and equi-edge mappings. Further discussion of the noncoordinate and coordinate bases can be found in \textcite{chen2024exocubed}. \par
The noncoordinate basis metric tensor is
\begin{equation}
    G_{ij} = 
    \begin{bmatrix}
        1  & \cos(\alpha) \\
        \cos(\alpha) & 1
    \end{bmatrix},
\label{eq:noncoord_metric_tensor}
\end{equation}

\noindent with a Jacobian determinant of $\sqrt{G} = \sqrt{\det(G_{ij})} = \sin(\alpha)$. The inverse metric tensor, $G^{ij} = (G_{ij})^{-1}$, is
\begin{equation}
    G^{ij} = 
    \frac{1}{\sin^2(\alpha)}
    \begin{bmatrix}
        1  & -\cos(\alpha) \\
        -\cos(\alpha) & 1
    \end{bmatrix}.
\end{equation}

As discussed in \textcite{nair2009diffusion}, the forms of the horizontal divergence, relative vorticity, and Laplacian operator acting on some scalar field $s$, in a curvilinear coordinate system of $(x^1,x^2)$, are given by 
\begin{equation}
    D = \frac{1}{\sqrt{G}}\bracked{\pfrac{}{x^1} \bracked{\sqrt{G} \tilde{u}} + \pfrac{}{x^2} \bracked{\sqrt{G} \tilde{v}} },
\end{equation}
\begin{equation}
    \zeta = \frac{1}{\sqrt{G}}\bracked{\pfrac{v}{x^1} - \pfrac{u}{x^2} },
\end{equation}
\begin{align}
\begin{split}
    \nabla^2 s = &\frac{1}{\sqrt{G}} \pfrac{}{x^1} \bracked{\sqrt{G} G^{11} \pfrac{s}{x^1} + \sqrt{G} G^{12} \pfrac{s}{x^2}} \\
    + &\frac{1}{\sqrt{G}} \pfrac{}{x^2} \bracked{\sqrt{G} G^{21} \pfrac{s}{x^1} + \sqrt{G} G^{22} \pfrac{s}{x^2}},
\end{split}
\end{align}

\noindent with ($u,v$) the covariant representation of the horizontal velocities and ($\tilde{u}, \tilde{v}$) the contravariant representation. \par
With the noncoordinate basis metric tensor (\ref{eq:noncoord_metric_tensor}) and cubed-sphere panel coordinates ($x,y$), the horizontal divergence is
\begin{equation}
    D = \frac{1}{\sin(\alpha)} \bracked{\pfrac{}{x}(\sin(\alpha) \tilde{u} ) + \pfrac{}{y}(\sin(\alpha) \tilde{v} )},
\label{eq:div_cubed_sphere}
\end{equation}

\noindent the relative vorticity is 
\begin{equation}
    \zeta = \frac{1}{\sin(\alpha)} \bracked{\pfrac{v}{x}- \pfrac{u}{y}},
\label{eq:vort_cubed_sphere}
\end{equation}

\noindent and the scalar Laplacian operator is

\begin{align}
\begin{split}
    \nabla^2 s = &\frac{1}{\sin(\alpha)} \pfrac{}{x} \left[ \frac{1}{\sin(\alpha)} \pfrac{s}{x} - \frac{\cos(\alpha)}{\sin(\alpha)} \pfrac{s}{y} \right] \\
    + &\frac{1}{\sin(\alpha)} \pfrac{}{y} \left[ - \frac{\cos(\alpha)}{\sin(\alpha)} \pfrac{s}{x} + \frac{1}{\sin(\alpha)} \pfrac{s}{y} \right].
\end{split}
\label{eq:full_Laplacian}
\end{align}

To discretise these expressions, we introduce central difference operators that are computed over one cell length,
\begin{subequations}
\begin{align}
    \delta_x (s_{i, j}) &= s_{i+0.5, j} - s_{i-0.5, j},  \\
    \delta_y (s_{i, j}) &= s_{i, j+0.5} - s_{i, j-0.5}.
\end{align}
\label{eq:central_diff_operator}
\end{subequations}

\noindent This allows the partial derivatives to be approximated by $\partial/\partial x \approx \delta_x/\Delta x, ~\partial/\partial y \approx \delta_y/\Delta y$, leading to
\begin{align}
\begin{split}
    \nabla^2 s = &\frac{1}{\Delta A} \delta_x \left[ \chi \frac{1}{\sin(\alpha)} \delta_x(s) - \frac{\cos(\alpha)}{\sin(\alpha)} \delta_y(s) \right] \\
    + &\frac{1}{\Delta A} \delta_y \left[ - \frac{\cos(\alpha)}{\sin(\alpha)} \delta_x(s) + \chi^{-1} \frac{1}{\sin(\alpha)} \delta_y(s) \right],
\end{split}
\label{eq:full_Laplacian_discrete}
\end{align}

\noindent with $\Delta A = \Delta x \Delta y \sin(\alpha)$ a first-order approximation for the cell area \parencite{santos2024analysis_phd}, and $\Delta x$ and $\Delta y$ approximated as constant over the Laplacian stencil. We have also introduced the cell aspect ratio of $\chi=\Delta y/\Delta x$ (Table \ref{table:cubed_sphere_key_quantities}). \par
A complication with using the discrete curvilinear Laplacian (\ref{eq:full_Laplacian_discrete}) for divergence and vorticity damping is the presence of cross-derivative terms, which increase the computational cost. Additionally, where the $\delta_x^2$ and $\delta_y^2$ operators take cell centre values and output to the cell centre, the $\delta_x\delta_y$ and $\delta_y \delta_x$ operators instead output to cell corners, which then requires interpolation to the cell centre. An alternative method to avoid the cross-derivative terms is to instead compute what we will term the \textit{pseudo-Laplacian}, as is used in FV3. Any terms pertaining to the pseudo-Laplacian will be denoted by a tilde. The pseudo-Laplacian is computed as  $\widetilde{\nabla}^2 s = \bm{\nabla} \cdot \bm{g}$, where $\bm{g} = \bm{\nabla} s = (\pfrac{s}{x},\pfrac{s}{y})$, so
\begin{equation}
    \widetilde{\nabla}^2 s = \frac{1}{\sin(\alpha)} \left[ \frac{\partial}{\partial x} \bracked{\sin(\alpha) \frac{\partial s}{\partial x} } + \frac{\partial}{\partial y} \bracked{\sin(\alpha) \frac{\partial s}{\partial y} } \right].
\label{eq:pseudo_Laplacian}
\end{equation}

\noindent Using the central difference operators (\ref{eq:central_diff_operator}) to approximate the gradient as $\bm{\nabla} s = (\frac{1}{\Delta x} \delta_x (s), \frac{1}{\Delta y} \delta_y (s))$ arrives at a discretisation for the pseudo-Laplacian,
\begin{equation}
    \widetilde{\nabla}^2 s = \frac{1}{\Delta A} \left[ \delta_x ( \chi \sin(\alpha)\delta_x (s)) +  \delta_y (\chi^{-1} \sin(\alpha) \delta_y (s)) \right].
\label{eq:pseudo_Laplacian_discrete}
\end{equation}

\noindent We further simplify this expression in preparation for the von Neumann stability analysis, which requires that the finite difference discretisation is linear and contains constant coefficients, as discussed in Chapter 3.2.2 of \textcite{Durran}. To freeze the coefficients, we approximate $\chi$ and $\sin(\alpha)$ as constant across the locations used in the central difference operators, leaving
\begin{equation}
    \widetilde{\nabla}^2 s = \frac{\sin(\alpha)}{\Delta A} \left[ \chi \delta_x^2 (s) +  \chi^{-1} \delta_y^2 (s) \right].
\label{eq:pseudo_Laplacian_discrete_simplified}
\end{equation}

\section{Linear stability analysis}
\label{section:linear_stab_analysis}
This section will begin with a linear stability analysis of a general order $2q$ damping operator that uses the pseudo-Laplacian, $\widetilde{\nabla}^{2q}s, ~q \in \mathbb{Z}^+$, before comparing its stability to that of the full Laplacian, $\nabla^{2q}s$. The analysis applies to both divergence and vorticity damping, using the dummy scalar $s \in \{D, \zeta \}$, with minor differences arising depending on the locations of $D$ and $\zeta$ on staggered Arakawa grids (Fig. \ref{fig:C_D_grids}). Note that as the local coordinates are the same on each cubed-sphere panel, an analysis of a single panel applies to the whole domain.

\subsection{Von Neumann analysis}
Consider a diffusion operator which uses the pseudo-Laplacian operator, with repeated applications for hyperviscosities,
\begin{align}
\begin{split}
    \widetilde{\nabla}^{2p} s = \frac{\sin(\alpha)}{\Delta A} \left[  \chi \delta^2_x \lbrac \widetilde\nabla^{2(p-1)} s \rbrac  + \chi^{-1} \delta^2_y \lbrac \widetilde\nabla^{2(p-1)} s  \rbrac  \right],
\label{eq:nabla_2q_s}
\end{split}
\end{align}

\noindent with $p \in \{1, .., q\}$. For example, sixth-order damping requires three applications of (\ref{eq:nabla_2q_s}) to compute $\widetilde\nabla^2s, \widetilde\nabla^4 s,$ then $\widetilde\nabla^6 s$. \par
Now, expand the central difference operators (\ref{eq:central_diff_operator}) in $\widetilde\nabla^{2q}s$ around an arbitrary primary grid index of $(i, j)$,
\begin{align}
\begin{split}
    \widetilde{\nabla}^{2q} s_{i,j} = &\frac{\sin(\alpha_{i,j})}{\Delta A_{i,j}} \\
    \Bigg{[}  &\chi_{i,j} [(\widetilde\nabla^{2(q-1)} s)_{i+1,j} - 2 (\widetilde\nabla^{2(q-1)} s)_{i,j} + (\widetilde\nabla^{2(q-1)} s)_{i-1,j}]\\
    + &\chi^{-1}_{i,j} [(\widetilde\nabla^{2(q-1)} s)_{i,j+1} - 2 (\widetilde\nabla^{2(q-1)} s)_{i,j} + (\widetilde\nabla^{2(q-1)} s)_{i,j-1}] \Bigg{]}.
\label{eq:deriv_step_A}
\end{split}
\end{align}

\noindent Hyperviscous $\widetilde{\nabla}^{2(q-1)}$ terms can be subsequently expanded with (\ref{eq:nabla_2q_s}) until no gradient operators remain on $s$. Each application of (\ref{eq:nabla_2q_s}) introduces an additional factor of $\sin(\alpha)/\Delta A$ and assumes that $\Delta x, \Delta y, \sin(\alpha)$ are constant across the stencil of the diffusion operator. \par
The expansion of $\widetilde{\nabla}^{2q} s$ (\ref{eq:deriv_step_A}) is now used in the time evolution of $s \in \{D,\zeta \}$ (\ref{eq:time_ev_eqs}), with a forward Euler discretisation of the time derivative,
\begin{align}
\begin{split}
    \frac{s^{n+1}_{i,j} - s^{n}_{i,j}}{\Delta t} = & 
    (-1)^{q+1} \nu_{s,2q} \frac{\sin(\alpha_{i,j})}{\Delta A_{i,j}} \\
    \Bigg{[}  &\chi_{i,j} [(\widetilde\nabla^{2(q-1)} s)^n_{i+1,j} - 2 (\widetilde\nabla^{2(q-1)} s)^n_{i,j} + (\widetilde\nabla^{2(q-1)} s)^n_{i-1,j}] \\
    + &\chi^{-1}_{i,j} [(\widetilde\nabla^{2(q-1)} s)^n_{i,j+1} - 2 (\widetilde\nabla^{2(q-1)} s)^n_{i,j} + (\widetilde\nabla^{2(q-1)} s)^n_{i,j-1}] \Bigg{]},
\end{split}
\end{align}

\noindent with the superscript $n$ or $n+1$ denoting the time index. Introducing the nondimensional damping coefficient (\ref{eq:nonD_damp_coeff}), gives
\begin{align}
\begin{split}
    s^{n+1}_{i,j} - s^{n}_{i,j} =& (-1)^{q+1} (C_{s, 2q} ~\Delta A_{\text{min}})^q \frac{\sin(\alpha_{i,j})}{\Delta A_{i,j}} \\
    \Bigg{[}  &\chi_{i,j} [(\widetilde\nabla^{2(q-1)} s)^n_{i+1,j} - 2 (\widetilde\nabla^{2(q-1)} s)^n_{i,j} + (\widetilde\nabla^{2(q-1)} s)^n_{i-1,j}] \\
    + &\chi^{-1}_{i,j} [(\widetilde\nabla^{2(q-1)} s)^n_{i,j+1} - 2 (\widetilde\nabla^{2(q-1)} s)^n_{i,j} + (\widetilde\nabla^{2(q-1)} s)^n_{i,j-1}] \Bigg{]}.
\end{split}
\label{eq:deriv_step_apply_coeff}
\end{align}

We now assume that the solution can be represented by Fourier modes, and express linear solutions to the divergence or vorticity time evolution equation (\ref{eq:time_ev_eqs}) at a point on the cubed-sphere as 
\begin{equation}
    s(x,y,t) = s_0 e^{\iota(k x + l y + \omega t)},
\label{eq:von_neumann}
\end{equation}

\noindent with $s_0$ the initial amplitude of the Fourier mode, $\mathbf{k}_H = (k,l)$ the horizontal wavenumber in the cubed-sphere panel coordinates, and $\iota=\sqrt{-1}$ denoting the imaginary unit. This allows the substitution of 
\begin{equation}
    s_{i,j}^n = s_0 e^{\iota(i k \Delta x + j l \Delta y + n \omega \Delta t )}
\end{equation}

\noindent into the finite difference discretisation (\ref{eq:deriv_step_apply_coeff}). Cancelling the common factor of $s_0 \exp (\iota(i k \Delta x + j l \Delta y + n \omega \Delta t))$ makes the square bracket term in (\ref{eq:deriv_step_apply_coeff}) become
\begin{equation}
    (-1)^q 4^q \frac{\sin^{q-1}(\alpha)}{(\Delta A)^{q-1}} \Bigg{[}  \chi \sin^2 \bracfrac{k \Delta x}{2} + \chi^{-1} \sin^2 \bracfrac{l \Delta y}{2} \Bigg{]}^q,
\end{equation}

\noindent The von Neumann representation of (\ref{eq:deriv_step_apply_coeff}) is then

\begin{align}
\begin{split}
    e^{\iota \omega \Delta t} - 1 =& \\ 
    - 4^q \lbrac \frac{C_{s, 2q} ~\Delta A_{\text{min}} \sin(\alpha) }{\Delta A} \rbrac^q \Bigg{[}  & \chi \sin^2 \bracfrac{k\Delta x}{2} + \chi^{-1} \sin^2 \bracfrac{l \Delta y}{2} \Bigg{]}^q.
\end{split}
\end{align}

To investigate numerical stability, we examine the growth rate of the diffusion. Letting $\widetilde{\Gamma}_{2q} = \exp(\iota \omega \Delta t)$ be the temporal amplification factor over one time step gives
\begin{equation}
    \widetilde\Gamma_{2q} (k \Delta x, l \Delta y) = 1 - \left[ \frac{4 C_{s, 2q} \Delta A_{\text{min}} \sin(\alpha)}{\Delta A} \lbrac \chi \sin^2 \bracfrac{k \Delta x}{2} + \chi^{-1} \sin^2 \bracfrac{l \Delta y}{2}  \rbrac \right]^q .
\label{eq:Gamma_2q_damp}
\end{equation}

\noindent We require that $|\widetilde{\Gamma}_{2q}| \leq 1$ to avoid exponential growth from the explicit diffusion; this is von Neumann stability. A similar statement for damping on uniform Cartesian meshes can be found in \textcite{klemp2017damping}, with the key differences in our cubed-sphere result (\ref{eq:Gamma_2q_damp}) being the inclusion of cell aspect ratios ($\chi$), the nonorthogonality term of $\sin(\alpha)$, and the nondimensional diffusion coefficient that uses the minimum cell area (\ref{eq:nonD_damp_coeff}). \par
The quantities of $k \Delta x$ and $l \Delta y$ in the linear stability expression (\ref{eq:Gamma_2q_damp}) are normalised wavenumbers in the cubed-sphere coordinates. Although a finer spatial resolution allows a larger range of wavenumbers to be resolved in a numerical model, the most oscillatory of these remains the 2$\Delta x$ wave at a normalised wavenumber of $(k\Delta x, l \Delta y) = (\pi, \pi)$. The $2 \Delta x$ wave changes sign between consecutive grid points and is the first mode to become unstable with excessive explicit diffusion. \par
The square bracketed term in the temporal amplification factor (\ref{eq:Gamma_2q_damp}) is positive-definite, so $\widetilde{\Gamma}_{2q} \leq 1, \forall q \geq 1$. Thus, linear instability only occurs when $\widetilde{\Gamma}_{2q} < -1$. The 2$\Delta x$ wave maximises the sine functions in (\ref{eq:Gamma_2q_damp}) to provide the greatest restriction on the damping,
\begin{equation}
    \widetilde{\Gamma}_{2q}(\pi, \pi) = 1 - \left[ \frac{4 C_{s, 2q} \Delta A_{\text{min}} \sin(\alpha)}{\Delta A} \lbrac \chi + \chi^{-1} \rbrac \right]^q .
\label{eq:2dx_amp_factor_full}
\end{equation}

To compare stability on the different gnomonic cubed-sphere grids, we combine the grid-specific parameters in (\ref{eq:2dx_amp_factor_full}) into a spatially-dependent \textit{grid stability function},
\begin{equation}
    \widetilde{\Psi}(x,y) = \frac{\Delta A}{\sin(\alpha) \Delta A_{\text{min}} (\chi + \chi^{-1})},
\label{eq:grid_stab_func_pseudo}
\end{equation}

\noindent which simplifies the 2$\Delta x$ wave temporal amplification factor (\ref{eq:2dx_amp_factor_full}) to
\begin{equation}
    \widetilde{\Gamma}_{2q}(\pi,\pi) = 1 - \lbrac \frac{4}{\widetilde\Psi} C_{s,2q}\rbrac^q.
\label{eq:2dx_amp_factor_simple}
\end{equation}

\noindent Enforcing von Neumann stability of $|\widetilde\Gamma_{2q}| \leq 1$ requires that
\begin{equation}
    C_{s,2q} \leq 2^{1/q} \frac{\widetilde\Psi_{\text{min}}}{4},
\label{eq:damp_stab_limit}
\end{equation}

\noindent where $\widetilde\Psi_{\text{min}}$ is the minimum value of the grid stability function over a cubed-sphere panel. $\widetilde{\Psi}_{\text{min}}$ is evaluated on the primary or offset grid according to the locations of the divergence or vorticity, which leads to small differences between explicit divergence and vorticity damping limits on the staggered Arakawa C- and D-grids. \par
A stronger restriction on the nondimensional damping coefficient can be made by ensuring that the amplification factor is non-negative. Although $\widetilde{\Gamma} \in [-1,0)$ is stable, it causes the wave amplitudes to change sign with each time step, which introduces unnecessary oscillations \parencite{klemp2017damping}. Choosing a coefficient of 
\begin{equation}
    C_{s,+} = \frac{\widetilde{\Psi}_{\text{min}}}{4}
\label{eq:stab_limit_osc_free}
\end{equation}

\noindent sets $\widetilde{\Gamma}(\pi,\pi) = 0 $ to ensure that all amplification factors are non-negative. We term (\ref{eq:stab_limit_osc_free}) the oscillation-free coefficient and note that its value is independent of the damping order, $q$. \par

\subsection{Grid comparisons}
\begin{figure}[htpb]
    \centering
     \includegraphics[width=\textwidth]{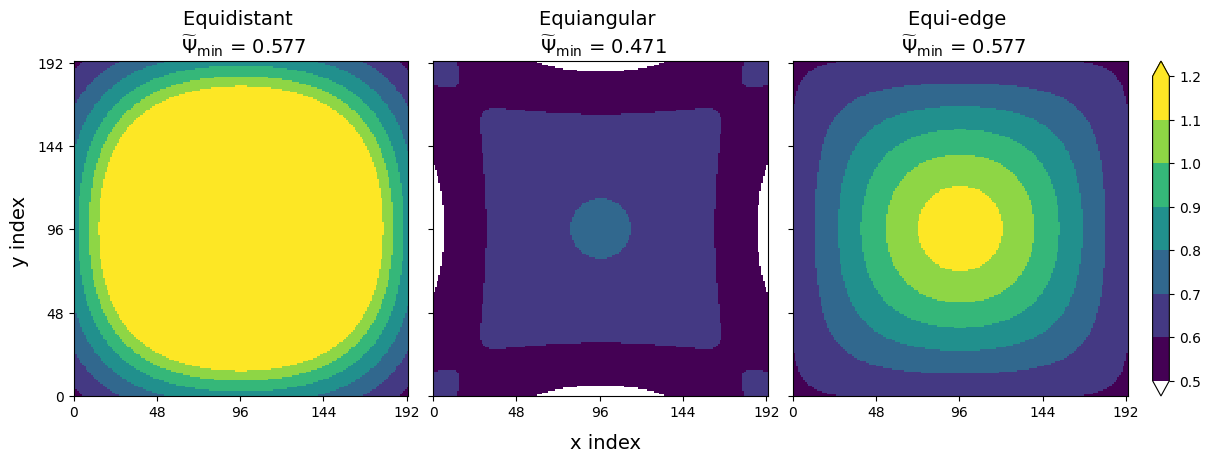}
    \caption[]{The grid stability function of $\widetilde\Psi$ evaluated on the offset grid for the three gnomonic cubed-sphere grids and a C192 resolution. We focus on its minimum value, $\widetilde\Psi_{\text{min}}$, which dictates the diffusive stability limit through (\ref{eq:damp_stab_limit}). A smaller $\widetilde\Psi_{\text{min}}$ for the equiangular grid (the lowest values are shown in white) corresponds to a stricter upper bound on the damping coefficient. $\widetilde\Psi_{\text{min}}$ is located at the smallest cells of each grid (Fig. \ref{fig:cubed_sphere_areas}); this is the middle of the panel edges for the equiangular grid, and the panel corners for the equidistant and equi-edge grids.}
    \label{fig:cubed_sphere_stab_func_div}
\end{figure}

We now compare analytical stability limits of the equidistant, equiangular and equi-edge mappings in Figure \ref{fig:cubed_sphere_stab_func_div}, by plotting their grid stability functions. For all mappings, $\widetilde{\Psi}_{\text{min}}$ occurs at the smallest cells, which are at the panel corners for the equidistant and equi-edge grids, and the middle of edges for the equiangular grid (Fig. \ref{fig:cubed_sphere_areas}). Thus, it is expected that linear instabilities will form at different points on the equiangular grid compared to the equidistant and equi-edge grids when using a pseudo-Laplacian damping operator. \par 
At the smallest cells, the grid stability function (\ref{eq:grid_stab_func_pseudo}) reduces to
\begin{equation}
    \widetilde{\Psi}_{\text{min}} = \frac{1}{\sin(\alpha) (\chi + \chi^{-1})}.
\label{eq:min_stab_func_C_grid}
\end{equation}

\noindent Further simplifications are possible for each cubed-sphere mapping. Consider the offset grid and an even number of cells per panel edge ($N_e$), which means that quantities are defined exactly at the panel corners and the middle of edges. At the middle of the panel edges on the equiangular grid, $\sin(\alpha) = 1$ and the aspect ratio is $\sqrt{2}$ or $1/\sqrt{2}$ \parencite{ronchi1996_equiangular}. This simplifies the grid stability function to
\begin{equation}
    \widetilde{\Psi}_{\text{min}} = \frac{\sqrt{2}}{3}.
\label{eq:stab_func_div_equi_ang}
\end{equation}

\noindent The equidistant and equi-edge grids have the smallest cells at the panel corners, where $\chi = 1$ and $\sin(\alpha) = \sqrt{3}/2$ (Table \ref{table:cubed_sphere_key_quantities}), so
\begin{equation}
    \widetilde{\Psi}_{\text{min}} = \frac{1}{\sqrt{3}}.
\label{eq:stab_func_div_equi_edge}
\end{equation}

\noindent Hence, the minimum value of the grid stability function is larger for the equidistant and equi-edge grids by a factor of $\sqrt{3/2} \approx 1.22$ compared to the equiangular grid. Accordingly, there are stricter linear stability limits on the equiangular grid for any order of explicit divergence or vorticity damping when using the pseudo-Laplacian operator. \par

On the primary grid (for $N_e$ even), quantities are offset from the corners or middle of panel edges by $\Delta x/2$ and/or $\Delta y/2$. This means that the simplified grid stability expressions (\ref{eq:stab_func_div_equi_ang}, \ref{eq:stab_func_div_equi_edge}) are not exact on the offset grid, but approximately hold. Table \ref{table:stab_func_primary_offset}, which compares $\widetilde\Psi_{\text{min}}$ on the primary and offset grids, shows that the differences are small and that $\widetilde\Psi_{\text{min}}$ on the primary grid tends to that of the offset grid with increasing resolution. Hence, the linear stability limits for explicit vorticity and divergence damping are extremely similar on the staggered C- and D-grids. \par

\begin{table}[htpb]
\caption[]{Comparing minimum evaluations of the grid stability function ($\widetilde\Psi_{\text{min}}$) on the primary and offset grids, to three decimal places. The offset grid minima are the same (for this level of accuracy) for the different horizontal resolutions shown. The minima for the primary grid are close to those on the offset grid, and there is increasing agreement with a finer resolution.}
\label{table:stab_func_primary_offset}
    \begin{center}
    \begin{tabular}{ |c||c|c|c||c|c|c|| }
    \hline
          & \multicolumn{3}{|c||}{Offset grid} & \multicolumn{3}{|c||}{Primary grid} \\
         \hline
          & C48 & C96 & C192 & C48 & C96 & C192 \\
         \hline
         \hline
        Equidistant & 0.577 & 0.577 & 0.577 & 0.573 & 0.575 & 0.576 \\
        \hline
        Equiangular & 0.471 & 0.471 & 0.471 & 0.474 & 0.473 & 0.472 \\
        \hline
        Equi-edge & 0.577 & 0.577 & 0.577  & 0.572 & 0.575 & 0.576 \\
        \hline
    \end{tabular}
    \end{center}
\end{table}

\begin{figure}
     \centering
    \includegraphics[width=0.7\textwidth]{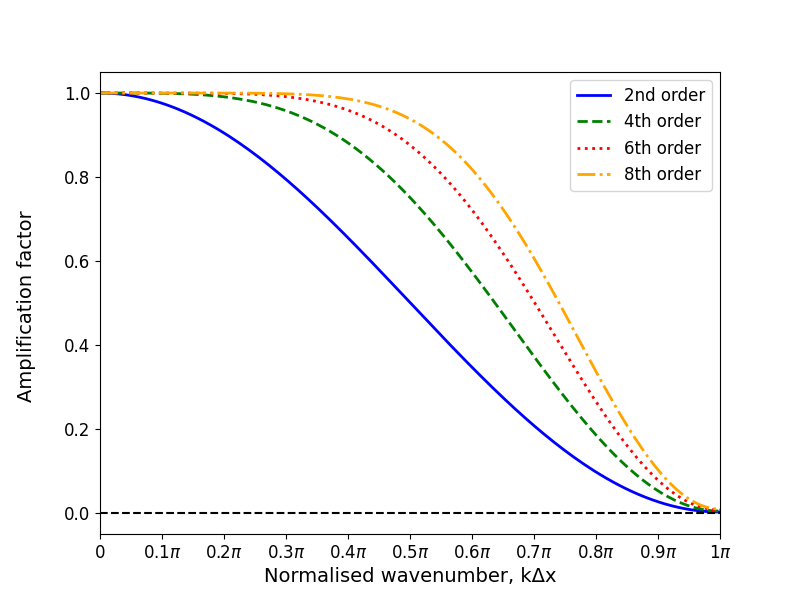}
        \caption{Amplification factors of $\widetilde \Gamma_{2q}$ against normalised wavenumber for the different orders of damping to show the impact of $q$ on scale selectivity. The oscillation-free stability coefficient (\ref{eq:stab_limit_osc_free}) and a C192 resolution are used. In contrast to the hyperviscosities, Laplacian diffusion significantly damps a wider range of wavenumbers. With increasing order, the damping becomes increasingly scale selective, with a smaller range of wavenumbers being noticeably diffused.}
    \label{fig:div_damp_strong_scale_select}
\end{figure}

We now examine amplification factors of $\widetilde\Gamma = \exp(\iota \omega \Delta t)$ at the smallest cell on each grid in Fig. \ref{fig:div_damp_strong_scale_select}. Amplification factors are shown along a diagonal slice of $k \Delta x = l \Delta y$ in normalised wavenumber space. This covers the range from the nonoscillatory $(k \Delta x, l \Delta y) = (0,0)$ wave, which is unaffected by explicit damping, to the most diffused $(k \Delta x, l \Delta y) = (\pi, \pi)$ wave. The oscillation-free coefficient (\ref{eq:stab_limit_osc_free}) is used, which (rounded down to 3 decimal places) is $C_{+} = 0.144$ for the equidistant and equi-edge grids, and $C_{+} = 0.117$ for the equiangular grid. Using oscillation-free coefficients means that the amplification factor for each wavenumber, in the smallest cell, is the same for all three gnomonic mappings. Fig. \ref{fig:div_damp_strong_scale_select} highlights that Laplacian diffusion damps most waves to an appreciable degree, including those of a low frequency. Increasing the order of diffusion for hyperviscosities shifts the amplification factor curve to the right, making the damping more scale-selective as more waves are minimally diffused.\par

\begin{figure}
    \centering
    \begin{subfigure}[b]{\textwidth}
        \includegraphics[width=\textwidth]{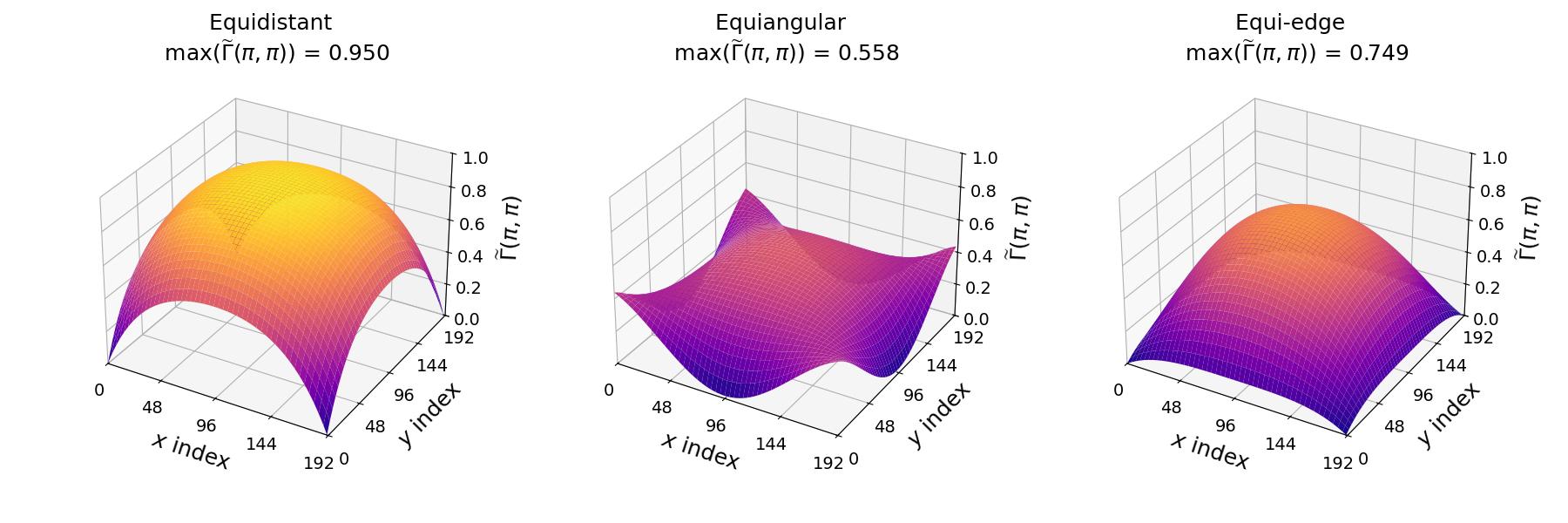}
        \caption{Fourth-order damping}
        \label{fig:2dx_damp_q2}
    \end{subfigure}
    \begin{subfigure}[b]{\textwidth}
        \includegraphics[width=\textwidth]{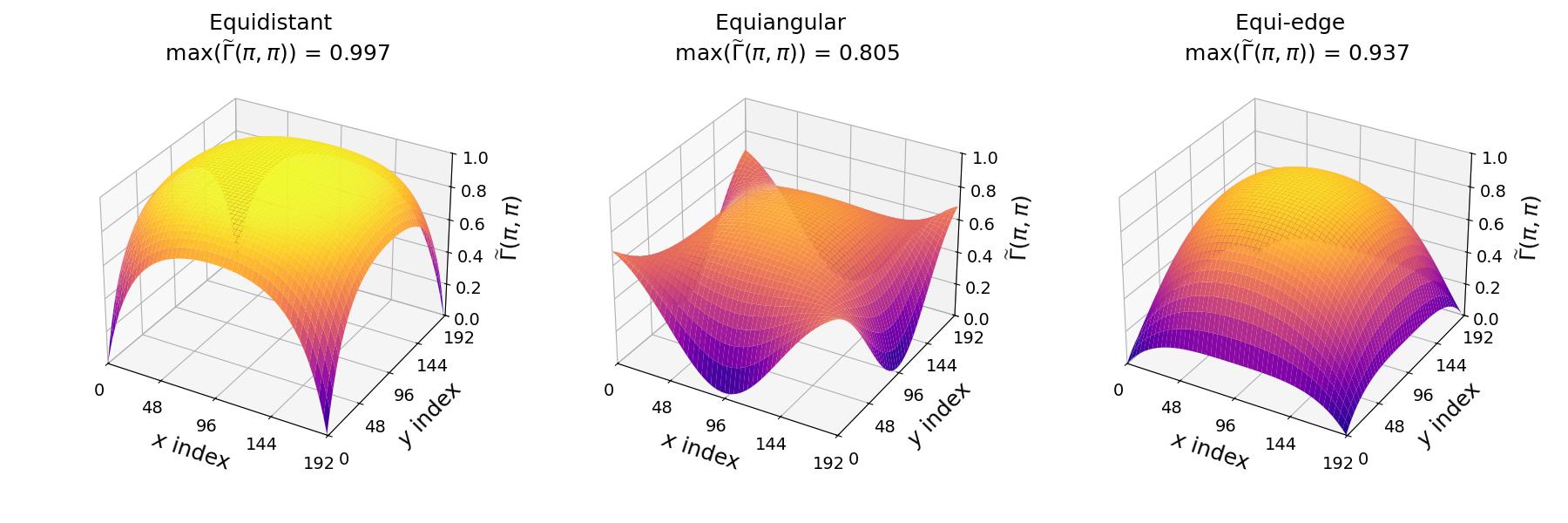}
        \caption{Eighth-order damping}
        \label{fig:2dx_damp_q4}
    \end{subfigure}
    \caption{Plots of the $2 \Delta x$ amplification factor, $\widetilde \Gamma(\pi, \pi)$, across a C192 panel of the cubed-sphere for the three gnomonic mappings. Oscillation-free coefficients are used. The amplification factor increases towards the centre of the panels due to increasing cell area and leads to less damping of the $2 \Delta x$ wave. Moving from fourth-order diffusion (a) to eighth-order diffusion (b) makes the damping weaker for all mappings.}
    \label{fig:2dx_damp_across_panel_pseudo}
\end{figure}

Next, we compare the damping of the $2 \Delta x$ wave across a cubed-sphere panel from the three gnomonic mappings. Figure \ref{fig:2dx_damp_across_panel_pseudo} shows the $2 \Delta x$ amplification factor of $\widetilde\Gamma(\pi,\pi)$ for fourth-order damping and eighth-order damping, again with the oscillation-free coefficient (\ref{eq:stab_limit_osc_free}). Amplification factors are always bounded above by one (\ref{eq:Gamma_2q_damp}), with larger values corresponding to less diffusion and $\widetilde\Gamma=1$ meaning no damping. $\widetilde\Gamma(\pi,\pi)$ is largest at the centre of the panels as a consequence of using the minimum cell area in the nondimensional diffusion coefficient (\ref{eq:nonD_damp_coeff}), as $\widetilde\Gamma(\pi,\pi)$ increases with a larger ratio of $\Delta A/\Delta A_{\text{min}}$ (\ref{eq:2dx_amp_factor_full}). With the largest range of cell areas, the equidistant grid has the weakest damping at the panel centre, whilst the equiangular grid has the most uniform damping of the 2$\Delta x$ wave across the panel. A higher order of damping increases $\widetilde\Gamma(\pi,\pi)$ and reduces the damping of the $2 \Delta x$ wave over the panel. Hence, a downside to increasing the order of diffusion, with the choice of nondimensional diffusion coefficient given by (\ref{eq:nonD_damp_coeff}), may be insufficient damping of the grid-scale waves in certain regions of the domain. This is especially true of the equidistant grid with eighth-order damping, where $\max(\widetilde\Gamma(\pi,\pi))=0.997$ (Fig. \ref{fig:2dx_damp_across_panel_pseudo}b).  \par

\clearpage

\subsection{Comparison between the full curvilinear Laplacian and pseudo-Laplacian operators}
As discussed in section \ref{sec:cubed_sphere_grids}c, the pseudo-Laplacian (\ref{eq:pseudo_Laplacian}) used in the previous von Neumann analysis is an approximation to the full curvilinear Laplacian (\ref{eq:full_Laplacian}). We now provide the stability restriction for the full Laplacian and compare it to that of the pseudo-Laplacian. For the cross-derivative terms of $\delta_x\delta_y$ and $\delta_y\delta_x$ (\ref{eq:full_Laplacian_discrete}), we take an equal average of evaluations at the four corners, leaving an expression entirely at the cell centre. The resulting amplification factor for the full Laplacian is
\begin{align}
\begin{split}
    \Gamma_{2q} (k \Delta x, l \Delta y) = 1 - \Bigg[ \frac{4 C_{s, 2q} \Delta A_{\text{min}}}{\Delta A \sin(\alpha)} &\Bigg( \chi \sin^2 \bracfrac{k \Delta x}{2} \\
    &-2 \cos(\alpha) \sin \bracfrac{k \Delta x}{2} \sin \bracfrac{l \Delta y}{2} \cos \bracfrac{k \Delta x}{2} \cos \bracfrac{l \Delta y}{2} \\
    &+ \chi^{-1} \sin^2 \bracfrac{l \Delta y}{2} \Bigg) \Bigg]^q .
\label{eq:Gamma_2q_damp_full_Laplacian}
\end{split}
\end{align}

\noindent The key differences from the amplification factor with the pseudo-Laplacian (\ref{eq:Gamma_2q_damp}) are the appearance of the cross-derivative term, which drops out when the grid is orthogonal and $\cos(\alpha) = 0$, and that the factorised metric term of $\sin(\alpha)$ now appears in the denominator instead of the numerator.  \par
If we again consider the $2 \Delta x$ wave, the cross-derivative term becomes zero, leaving
\begin{equation}
    \Gamma_{2q}(\pi, \pi) = 1 - \left[ \frac{4 C_{s, 2q} \Delta A_{\text{min}}}{\Delta A \sin(\alpha)} \lbrac \chi + \chi^{-1} \rbrac \right]^q.
\label{eq:2dx_amp_factor_full_Laplacian}
\end{equation}

\noindent The corresponding grid stability function is
\begin{equation}
    \Psi(x,y) = \frac{\Delta A \sin(\alpha)}{ \Delta A_{\text{min}} (\chi + \chi^{-1})}.
\label{eq:grid_stab_func_full}
\end{equation}

\noindent Hence, the relationship between the grid stability function from the full Laplacian (\ref{eq:grid_stab_func_full}) and pseudo-Laplacian (\ref{eq:grid_stab_func_pseudo}) is
\begin{equation}
    \Psi(x,y) = \sin^2(\alpha) \widetilde{\Psi} (x,y).
\label{eq:grid_stab_relation}
\end{equation}

\noindent This means that the grid stability functions are equal where the grid is orthogonal and $\Psi < \widetilde \Psi$ elsewhere. \par
\begin{figure}
    \centering
    \includegraphics[width=\linewidth]{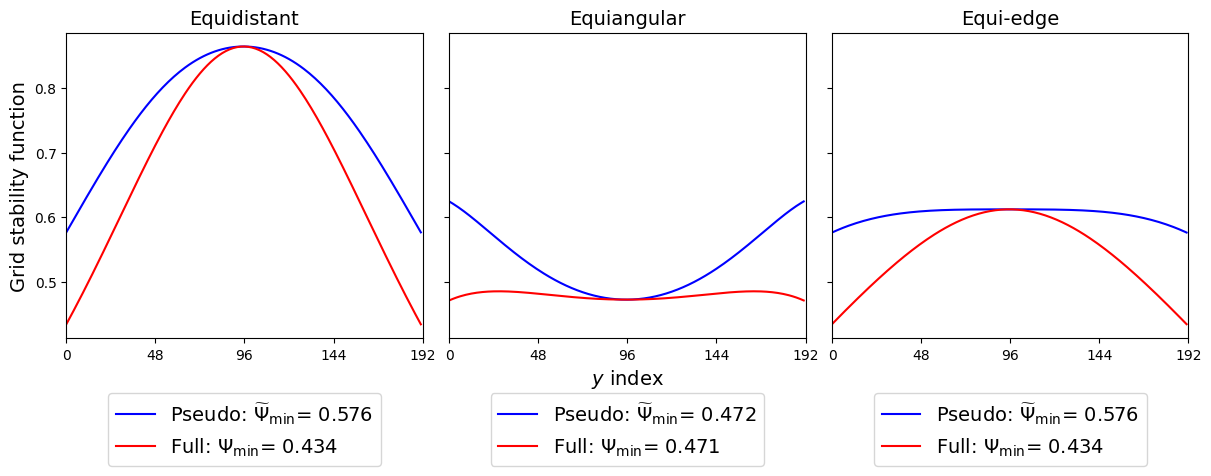}
    \caption{Comparing grid stability functions from the pseudo- and full Laplacians along an edge of a cubed-sphere panel, which includes the important locations of the corners and middle of a panel edge. A C192 resolution is used. Using the full Laplacian significantly reduces the grid stability function for the equidistant and equi-edge grids, with the largest reduction at the cell corner. The reduction of $\Psi$ at the cell corners for the equiangular grid means that the middle of the panel edges are no longer the most constraining locations, although $\Psi_{\text{min}} \approx \widetilde \Psi_{\text{min}}$.}
    \label{fig:compare_grid_stab_panel_edge}
\end{figure}

Figure \ref{fig:compare_grid_stab_panel_edge} compares the grid stability functions along an arbitrary panel edge. At the middle of a panel edge, the grid is orthogonal and $\Psi=\widetilde \Psi$, whereas at the panel corners, $\Psi= 0.75 \widetilde \Psi$. This reduction of the grid stability function when moving to the full Laplacian leads to $\Psi_{\text{min}}$ being at the panel corners of the equiangular grid, in contrast to $\widetilde \Psi_{\text{min}}$ at the middle of panel edges (Fig \ref{fig:cubed_sphere_stab_func_div}). For the equi-edge and equidistant grids, where $\widetilde \Psi_{\text{min}}$ is already at the panel corners, the grid stability function is reduced by 25\% when using the full Laplacian instead of the pseudo-Laplacian. Additionally, where $\widetilde \Psi_{\text{min}}$ is smallest on the equiangular grid, $\Psi_{\text{min}}$ is smallest on the equidistant and equi-edge grids, which means that smaller damping coefficients are required for stability (\ref{eq:damp_stab_limit}). \par

\begin{figure}
    \centering
    \begin{subfigure}[b]{\textwidth}
        \includegraphics[width=\textwidth]{two_dx_damp_over_panel_q2_pseudo.jpg}
        \caption{Pseudo-Laplacian operator}
        \label{fig:2dx_damp_pseudo}
    \end{subfigure}
    \begin{subfigure}[b]{\textwidth}
        \includegraphics[width=\textwidth]{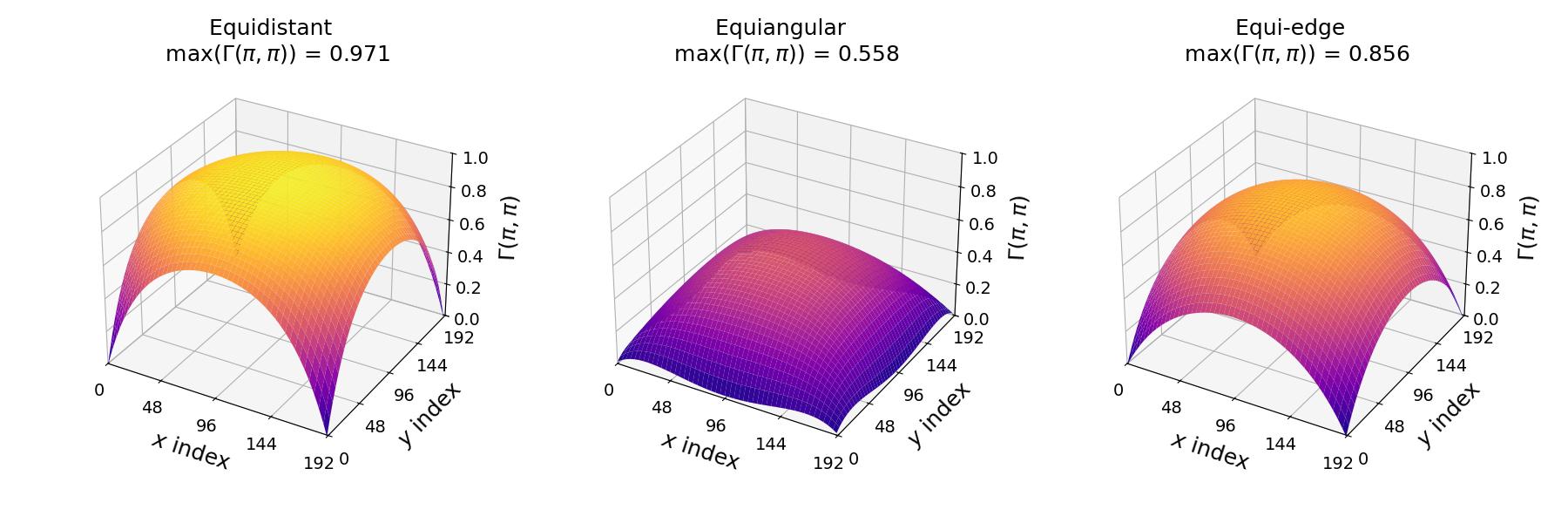}
        \caption{Full Laplacian operator}
        \label{fig:2dx_damp_full}
    \end{subfigure}
    \caption{Amplification factors of the $2 \Delta x$ wave across a C192 panel with fourth-order divergence damping, comparing the pseudo- and full Laplacian operators ((a) is the same as Fig. \ref{fig:2dx_damp_across_panel_pseudo}(a)). Oscillation-free coefficients (\ref{eq:stab_limit_osc_free}) are used, with $C_+ = \Psi_{\text{min}}/4$ smaller for the full Laplacian operator compared to $C_+ = \widetilde \Psi_{\text{min}}/4$ for the pseudo-Laplacian. The use of the full Laplacian increases the damping inhomogeneity of the $2 \Delta x$ wave for the equidistant and equiangular grids. Whilst the spatial layout of amplification factors changes for the equiangular grid, the range of the damping strength is unchanged.}
    \label{fig:comp_2dx_pseudo_full}
\end{figure}

We now compare how the full Laplacian impacts the damping of the $2 \Delta x$ wave across a cubed-sphere panel, shown for fourth-order damping in Fig. \ref{fig:comp_2dx_pseudo_full}. The reduction in $\Psi$ where the panel is nonorthogonal (\ref{eq:grid_stab_relation}) leads to greater inhomogeneity in the $2 \Delta x$ damping on the equidistant and equi-edge grids, with a larger range of amplification factors with the full Laplacian compared to the pseudo-Laplacian. By contrast, the range of amplification factors for the equiangular grid is virtually the same with both operators. \par

To summarise, the use of the full Laplacian instead of the pseudo-Laplacian requires smaller damping coefficients for stability on the equidistant and equi-edge grids. Whilst it has minimal impact on the stability limit for the equiangular grid, it shifts the predicted location of instability from the middle of panel edges to the panel corners. The use of the full Laplacian also increases the inhomogeneity in damping the $2 \Delta x$ wave on the equi-edge grid (Fig. \ref{fig:comp_2dx_pseudo_full}), which may be an additional motivation for the use of the pseudo-Laplacian in FV3. In the next sections, we will focus on the pseudo-Laplacian, as we compare the theoretical stability limits to practical damping restrictions in the hydrostatic CAM-FV3 dynamical core, using the equiangular and equi-edge grids.

\clearpage

\section{Application to CAM-FV3}
\label{sec:fv3_and_test_case}

\subsection{FV3 overview}
FV3 primarily uses the Arakawa D-grid to prioritise accurate rotational dynamics, whilst performing transport on the C-grid with interpolated winds for advection. Hence, we use the offset grid stability limits for divergence damping and the primary grid for vorticity damping (Fig. \ref{fig:C_D_grids}). FV3 uses the finite volume transport scheme of \textcite{lin1996multidimensional} and a floating vertical Lagrangian coordinate \parencite{lin2004vertically}. This means that although we will be discussing horizontal momentum equations, these are not technically horizontal, but apply along the two-dimensional surfaces prescribed by each Lagrangian level. \par
FV3's horizontal momentum equations are defined in Equations (6.1d) and (6.1e) in \textcite{harris2021scientific_FV3_doc} as
\begin{subequations}
\begin{align}
    \pfrac{u}{t} &= \frac{u^{n+1}-u^n}{\Delta t} = (\mathcal{T}_Y +  \mathcal{V}_{y,2q}) - \delta_x (K^* - \mathcal{D}_x) + P_x, \\
    \pfrac{v}{t} &= \frac{v^{n+1}-v^n}{\Delta t} = -(\mathcal{T}_X + \mathcal{V}_{x,2q}) - \delta_y (K^* - \mathcal{D}_y) + P_y,
\end{align}
\label{eq:fv3_momentum_eqs}
\end{subequations}

\noindent where $\Delta t$ is the dynamics time step. The first bracketed terms on the right-hand side pertain to two-dimensional transport on each Lagrangian surface: $\mathcal{T}_X, \mathcal{T}_Y$ are the transport operators of \textcite{lin1996multidimensional} which contain implicit diffusion, whilst $\mathcal{V}_{x,2q}, \mathcal{V}_{y,2q}$ are diffusive fluxes that implement the optional vorticity damping. In the next term, $K^*$ is the kinetic energy, and $\mathcal{D}_x, \mathcal{D}_y$ represent explicit divergence damping. Lastly, $P_x, P_y$ denote horizontal pressure gradients. In the hydrostatic model, all terms in the horizontal momentum equations are evaluated explicitly using the time step $n$ fields, meaning that the same time step size is used for diffusion and transport. \par

\subsection{Test cases}
\label{subsec:test_cases}
To examine the stability of divergence and vorticity damping in CAM-FV3, we identify stable coefficients in two test cases. Although simpler than a full weather or climate simulation, these tests are sufficient for investigating diffusive instabilities that will be present in more complex flows. The CAM-FV3 default settings are mostly untouched, with our focus on modifying diffusion parameters. To find the practical damping limit, the largest value of $C_{D, 2q}$ or $C_{\zeta, 2q}$ (to three decimal places) is sought such that the simulation completes without numerical blow-up. \par
The first test case is the idealised baroclinic wave of \textcite{jablonowski2006baroclinic}, referred to as the JW2006 test. The initial condition prescribes a background steady-state that is a solution to the adiabatic governing equations. An added perturbation then triggers a baroclinic instability, which is an important mechanism for generating mid-latitude weather patterns. We run this test for fifteen days at resolutions of C96, C192, and C384, with dynamics time steps of $\Delta t = 150 ~\text{s}, ~\Delta t = 75 ~\text{s},$ and $\Delta t = 18.75 ~\text{s}$, respectively; the reduction of $\Delta t$ by a factor of four at the finest resolution is required for general model stability. \par
The second test case is a spun-up Held-Suarez simulation \parencite{held1994proposal}, referred to as the HS1994 test. The Held-Suarez case is commonly used in climate applications to study general atmospheric circulation. The spin-up simulation is run for 360 days, with the final state being stored and used as the initial condition for the diffusion tests, which are run for another thirty days. The spin-up process leads to a more complex flow field than in the JW2006 test. The diffusion tests use a C96 resolution and a dynamics time step of $\Delta t = 150 ~\text{s}$. The Held-Suarez relaxation terms, implemented as physics forcings, are applied every twelve dynamics time steps. \par
The namelist parameter of \texttt{grid\_type} selects the gnomonic cubed-sphere mapping in CAM-FV3. We compare the equiangular and equi-edge grids, but not the equidistant grid, due to it being unsupported in the CAM version used. We will examine five combinations of horizontal transport schemes, which covers the available choices in CAM6. There are two unlimited schemes --- \textit{virtually-inviscid} and \textit{intermediate unlimited} --- and three monotonic schemes --- \textit{CAM hydrostatic default}, \textit{Lin monotonic}, and \textit{Huynh monotonic}. The appendix provides further detail on these transport schemes. \par

\subsection{Divergence damping}
From \textcite{harris2021scientific_FV3_doc}, the pseudo-Laplacian or corresponding hyperviscous operator acting on horizontal divergence in FV3 is
\begin{align}
\begin{split}
    \widetilde \nabla^{2q} D = \frac{1}{\Delta A_c} \bigg[ &\delta_x \lbrac \frac{\Delta y_c}{\Delta x_d} \sin(\alpha) \delta_x \bracked{\widetilde\nabla^{2(q-1)} D } \rbrac \\
    + &\delta_y \lbrac \frac{\Delta x_c}{\Delta y_d} \sin(\alpha) \delta_y \bracked{\widetilde\nabla^{2(q-1)} D} \rbrac   \bigg].\label{eq:FV3_div_laplace_and_higher_operator}
    \end{split}
\end{align}

\noindent Note the presence of both primary ($_d$) and offset ($_c)$ grid quantities. To apply the linear stability limit (\ref{eq:damp_stab_limit}) we assume that $\Delta x_c = \Delta x_d = \Delta x$ and $\Delta y_c = \Delta y_d = \Delta y$, which is in line with the approximation in Section \ref{section:linear_stab_analysis} that $\Delta x$ and $\Delta y$ are constant over the diffusive stencil. \par
FV3 allows for the divergence damping terms ($\mathcal{D}_x,\mathcal{D}_y$) to be solely Laplacian ($q=1$) or hyperdiffusive ($q \in \{ 2, 3, 4 \}$), or a combination of these,
\begin{subequations}
    \begin{align}
    \mathcal{D}_x &= \frac{\nu_{D,2}}{\Delta x_d} D + (-1)^{q+1} \frac{\nu_{D,2q}}{\Delta x_d} \widetilde\nabla^{2(q-1)}  D, \\
    \mathcal{D}_y &= \frac{\nu_{D,2}}{\Delta y_d} D + (-1)^{q+1} \frac{\nu_{D,2q}}{\Delta y_d} \widetilde\nabla^{2(q-1)}  D,
\end{align}
\label{eq:fv3_curly_D_defs}
\end{subequations}

\noindent which in the FV3 horizontal momentum equations (\ref{eq:fv3_momentum_eqs}) leads to a time-evolution equation for the divergence of 
\begin{equation}
    \pfrac{D}{t} = ... + \nu_{D,2} \widetilde\nabla^2 D + (-1)^{q+1} \nu_{D,2q} \widetilde\nabla^{2q} D.
    \label{eq:FV3_D_time_ev}
\end{equation}

\noindent Pure hyperdiffusion is obtained by zeroing the Laplacian coefficient of $\nu_{D,2}$ and pure Laplacian diffusion occurs when $\nu_{D,2q} = 0$. A stability limit for mixed Laplacian and hyperviscous divergence damping combines $2 \Delta x$ wave amplification factors (\ref{eq:2dx_amp_factor_simple}) from both diffusive terms,
\begin{equation}
    \widetilde\Gamma_{D, 2+2q}(\pi, \pi) = 1 - \frac{4}{\widetilde\Psi_{c,\text{min}}} C_{D,2} - \lbrac \frac{4}{\widetilde\Psi_{c,\text{min}}} C_{D,2q} \rbrac^q,
\end{equation}

\noindent with $\widetilde\Psi_{c,\text{min}}$ denoting the minimum value of the pseudo-Laplacian grid stability function (\ref{eq:grid_stab_func_pseudo}) on the offset C-grid, and $C_{D,2}, C_{D,2q}$ the nondimensional diffusion coefficients (\ref{eq:nonD_damp_coeff}) used in CAM-FV3. Von Neumann stability requires that
\begin{equation}
    \frac{4}{\widetilde\Psi_{c,\text{min}}} C_{D,2} + \lbrac \frac{4}{\widetilde\Psi_{c,\text{min}}} C_{D,2q} \rbrac^q \leq 2.
\label{eq:stab_mixed_div_damp}
\end{equation}

\noindent This means that selecting the magnitude of either the Laplacian or hyperviscous coefficient constrains the other for linear stability. \par
\begin{figure}
     \centering
     \begin{subfigure}[b]{0.49\textwidth}
         \centering
         \includegraphics[width=19pc]{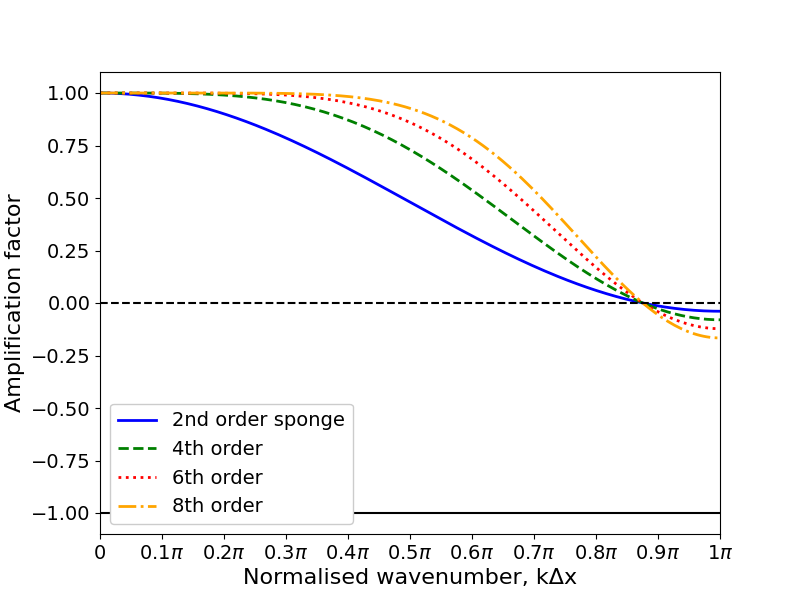}
         \caption{Equi-edge}
         \label{fig:equi_edge_scale_select}
     \end{subfigure}
     \hfill
     \begin{subfigure}[b]{0.49\textwidth}
         \centering
         \includegraphics[width=19pc]{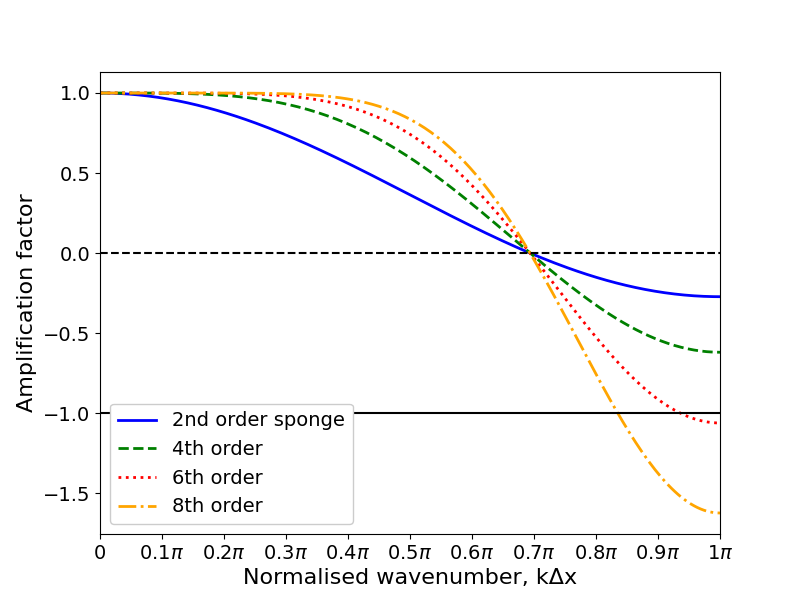}
         \caption{Equiangular}
         \label{fig:equi_ang_scale_select}
     \end{subfigure}
        \caption{Amplification factors against normalised wavenumber on the equi-edge and equiangular grids at a C192 resolution. Coefficients of $C_{D, 2q} = 0.15$ and $C_{D, 2}  = 0.15$ are used, which are the defaults for the hyperviscosity and the strongest Laplacian sponge layer, respectively. All CAM-FV3 defaults are stable for the equi-edge grid, whereas these choices lead to linear instability on the equiangular grid for sixth- and eighth-order divergence damping, as the amplification factor drops below -1.}
    \label{fig:CAM_def_scale_select}
\end{figure}

We now examine divergence damping amplification factors with the CAM-FV3 default coefficients. In the majority of the computational domain, the default is solely hyperviscous divergence damping, with $C_{D,2}=0, C_{D,2q} = 0.15$. FV3 also contains two sponge layers at the model top, which only use Laplacian diffusion. The uppermost sponge layer (covering the top two model levels) is the strongest, with $C_{D,2}=0.15$. Hence, the amplification factors in Fig. \ref{fig:CAM_def_scale_select} use $C_{D,2}=0.15$ for Laplacian diffusion and $C_{D,2q}=0.15$ for fourth-, sixth-, and eighth-order hyperviscosities. FV3 also uses Rayleigh damping in the upper model levels \parencite{harris2021scientific_FV3_doc}, but we will not consider this here. \par
The CAM-FV3 default damping coefficients are stable for all orders of divergence damping on the equi-edge grid (Fig. \ref{fig:CAM_def_scale_select}a). However, there are negative amplification factors for the largest wavenumbers as the coefficients are above the oscillation-free value of $C_{D,+} = 0.144$. The equiangular grid with CAM-FV3 default coefficients is stable for Laplacian and fourth-order divergence damping, but is unstable for sixth- and eighth-order (Fig. \ref{fig:CAM_def_scale_select}b). As CAM-FV3 uses eighth-order divergence damping by default, switching to the equiangular grid without reducing the divergence damping strength will likely instigate a simulation blow-up. Accordingly, the Held-Suarez spin-up simulation on the equiangular grid used a reduced $C_{D,2q} = 0.13$. \par

We now move to the CAM-FV3 simulations, with Table \ref{table:div_damp_limits} quoting the largest stable divergence damping coefficients obtained in the JW2006 and HS1994 tests on the equi-edge and equiangular grids. There is good agreement between linear theory and the practical limits in both tests and across the resolutions. In many cases, the practical limit is slightly above the analytical prediction; these simulations may still contain a linear instability that would cause numerical blow-up at longer times. The closeness of the practical limits between the JW2006 and HS1994 tests indicates that the linear theory should generalise to an arbitrary atmospheric flow. \par
As expected from the von Neumann analysis of the pseudo-Laplacian operator, the upper bounds on $C_{D,2q}$ are larger for the equi-edge grid than the equiangular grid. Stable values for the equi-edge grid are in line with the advice in \textcite{harris2021scientific_FV3_doc} that $C_{D, 2q} \leq 0.16$ should ensure numerical stability. The upper bounds for some sixth-order and all eighth-order tests on the equiangular grid are below 0.15, reiterating the potential for instability with the CAM-FV3 default of $C_{D, 2q} = 0.15$ and the equiangular gnomonic mapping (Fig. \ref{fig:equi_ang_scale_select}). Increasing the horizontal resolution in the JW2006 test from C96 to C192 or C384 leads to a closer agreement with linear theory. Accordingly, we expect that the analytical stability limits will remain accurate at the finer horizontal resolutions that are used for numerical weather prediction. \par 
The divergence damping limits are impacted by the choice of horizontal transport scheme, with the two unlimited methods tending to have a stricter upper bound (Table \ref{table:div_damp_limits}). This indicates that although the divergent modes are theoretically `invisible' to the FV3 transport scheme \parencite{harris2021scientific_FV3_doc}, this may not be entirely the case in practice. For eighth-order damping with virtually-inviscid transport, no tested value of $C_{D, 2q}$ enabled a stable simulation at C384 on both grids and C192 on the equiangular grid --- most likely, the amount of explicit damping required to stabilise the divergent modes exceeds the linear stability limit. Additionally, no stable $C_{D, 2q}, q \in \{ 2,3 \}$ could be found with the Huynh transport scheme in the Held-Suarez test, and the reason for this is unclear.  \par

\begin{sidewaystable}[htpb]
\caption[]{The largest values of $C_{D,2}$ and $C_{D,2q}$ that can be used for CAM-FV3 divergence damping in the JW2006 and HS1994 tests, to three decimal places. There is good agreement with the linear stability limits, given in the first row, which were computed using (\ref{eq:damp_stab_limit}) and rounded down to three decimal places. Differences in the upper bounds arise from the choice of the five horizontal transport schemes. Configurations where no stable divergence damping coefficient could be found are denoted by `None'. }
\label{table:div_damp_limits}
    \begin{center}
    \begin{tabular}{ |c|c|c||c|c|c|c||c|c|c|c| }
    \hline
          &  & & \multicolumn{4}{|c|}{Equi-edge} & \multicolumn{4}{|c|}{Equiangular} \\
         \hline
         \hline
          &  & Order of divergence damping & 2nd & 4th & 6th & 8th & 2nd & 4th & 6th & 8th \\
         \hline
        \hline
        Resolution & Test case & Linear stability limit & 0.288 & 0.204 & 0.181 & 0.171 & 0.235 & 0.166 & 0.148 & 0.140\\
        \hline
        \hline
        C96 & JW2006 & Default CAM (monotonic) & 0.295 & 0.206 & 0.184 & 0.174 & 0.242 & 0.171 & 0.153 & 0.144 \\
        \hline
        C96 & JW2006 & Lin monotonic & 0.295 & 0.206 & 0.184 & 0.174 & 0.241 & 0.171 & 0.152 & 0.144 \\
        \hline
        C96 & JW2006 & Huynh monotonic & 0.294 & 0.206 & 0.184 & 0.174 & 0.242 & 0.171 & 0.153 & 0.144 \\
        \hline
        C96 & JW2006 & Virtually-inviscid unlimited & 0.296 & 0.203 & 0.183 & 0.173 & 0.241 & 0.169 & 0.151 & 0.143 \\
        \hline
        C96 & JW2006 & Intermediate unlimited & 0.291 & 0.203 & 0.182 & 0.173 & 0.242 & 0.171 & 0.152 & 0.144 \\
        \hline
        \hline
        C96 & HS1994 & Default CAM (monotonic) & 0.296 & 0.205 & 0.184 & 0.174 & 0.243 & 0.170 & 0.152 & 0.143 \\
        \hline
        C96 & HS1994 & Lin monotonic & 0.296 & 0.205 & 0.183 & 0.173 & 0.243 & 0.169 & 0.152 & 0.143 \\
        \hline
        C96 & HS1994 & Huynh monotonic & 0.296 & 0.205 & None & None & 0.243 & 0.170 & None & None \\
        \hline
        C96 & HS1994 & Virtually-inviscid unlimited & 0.292 & 0.200 & 0.181 & 0.172 & 0.241 & 0.166 & 0.150 & 0.142 \\
        \hline
        C96 & HS1994 & Intermediate unlimited & 0.295 & 0.201 & 0.181 & 0.172 & 0.240 & 0.166 & 0.150 & 0.142 \\
        \hline
        \hline
        \hline
        C192 & JW2006 & Default CAM (monotonic) & 0.291 & 0.204 & 0.182 & 0.172 & 0.239 & 0.169 & 0.151 & 0.142 \\
        \hline
        C192 & JW2006 & Lin monotonic & 0.291 & 0.204 & 0.182 & 0.172 & 0.238 & 0.169 & 0.151 & 0.142 \\
        \hline
        C192 & JW2006 & Huynh monotonic & 0.291 & 0.204 & 0.182 & 0.172 & 0.239 & 0.169 & 0.151 & 0.142 \\
        \hline
        C192 & JW2006 & Virtually-inviscid unlimited & 0.285 & 0.200 & 0.180 & 0.170 & 0.234 & 0.167 & 0.149 & None \\
        \hline
        C192 & JW2006 & Intermediate unlimited & 0.284 & 0.201 & 0.180 & 0.171 & 0.239 & 0.168 & 0.150 & 0.142\\
        \hline
        \hline
        C384& JW2006 & Default CAM (monotonic) & 0.290 & 0.205 & 0.183 & 0.173 & 0.238 & 0.168 & 0.150 & 0.142 \\
        \hline
        C384 & JW2006 & Virtually-inviscid unlimited & 0.289 & 0.203 & 0.182 & None & 0.235 & 0.167 & 0.149 & None \\
        \hline
        \hline
        C96 & JW2006 & Default CAM, Flow-dependent & N/A & 0.205 & 0.184 & 0.174 & N/A & 0.168 & 0.151 & 0.143 \\
        \hline
        C96 & JW2006 & Virtually-inviscid, Flow-dependent & N/A & 0.200 & 0.181 & 0.172 & N/A & 0.165 & 0.148 & 0.140 \\
        \hline
    \end{tabular}
    \end{center}
\end{sidewaystable}

\begin{figure}[htpb]
    \centering
    \begin{subfigure}[b]{\textwidth}
        \includegraphics[width=\textwidth]{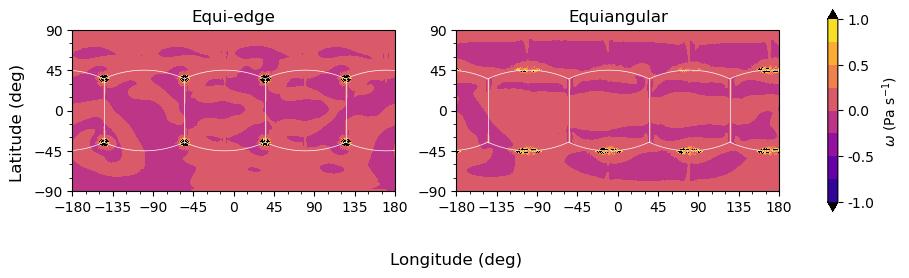}
        \caption{Locations of diffusive instability in the JW2006 test, at 850 hPa. These plots show fields at 0.63 days for the equi-edge grid and 0.79 days for the equiangular grid.}
        \label{fig:blowup_loc_jw2006}
    \end{subfigure}
    \begin{subfigure}[b]{\textwidth}
        \includegraphics[width=\textwidth]{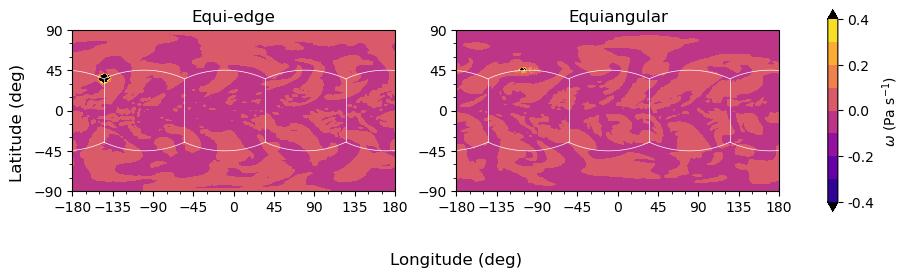}
        \caption{Locations of diffusive instability in the HS1994 test, at the lowest model level with reference pressure of $\approx$993 hPa. These plots show fields at 0.85 days for the equi-edge grid and 0.52 days for the equiangular grid.}
        \label{fig:blowup_loc_hs1994}
    \end{subfigure}
    \caption[]{The vertical pressure velocity ($\omega$) field at the last computable output time step in diffusively-unstable CAM-FV3 simulations. These use the default horizontal transport scheme, a C96 resolution, and sixth-order divergence damping, with $C_{D,6}$ set as 0.001 greater than the maximum stable value quoted in Table \ref{table:div_damp_limits}. The largest values of $|\omega|$ are shown in black to highlight the location(s) of numerical instability on the equi-edge and equiangular grids. The cubed-sphere panel edges are shown in white. FV3 aligns the panel edges at longitudinal values of $-145^{\circ}, -55^{\circ}, 35^{\circ}, 125^{\circ}$, which is a $10^{\circ}$ shift from symmetry about $0^{\circ}$ to avoid placing a panel edge over the mountains of Japan. In both tests, we see the expected locations of instability with damping operators that use the pseudo-Laplacian: panel corners for the equi-edge grid, and the middle of panel edges for the equiangular grid.}
    \label{fig:grid_blowup_comp}
\end{figure}

We next examine the locations of blow-up on the equi-edge and equiangular grids. Figure \ref{fig:grid_blowup_comp} shows the vertical pressure velocity field at the last completed output time step from unstable JW2006 and HS1994 simulations. Sixth-order divergence damping is used, with $C_{D,6}$ set at 0.001 greater than the maximum stable value identified in Table \ref{table:div_damp_limits} for the C96 resolution. The regions of numerical instability are seen to be the panel corners on the equi-edge grid and the middle of panel edges on the equiangular grid. These are the locations of the smallest cells and $\widetilde \Psi_{\text{min}}$ (Fig. \ref{fig:cubed_sphere_stab_func_div}), confirming the expected locations of instability from the von Neumann analysis of the pseudo-Laplacian. In the JW2006 test, the instability grows more uniformly at the smallest cells across the grid, as clearly seen for all eight corners on the equi-edge grid (Fig. \ref{fig:grid_blowup_comp}a), whereas for the HS1994 test, the diffusive instability grows most strongly at a single region of the cubed-sphere (Fig. \ref{fig:grid_blowup_comp}b).

Next, we investigate the simultaneous application of Laplacian and hyperviscous divergence damping, where the superposition of amplification factors needs to be bounded for von Neumann stability (\ref{eq:stab_mixed_div_damp}). For a fixed Laplacian coefficient, the linear stability restriction on the hyperviscosity becomes
\begin{equation}
    C_{D,2q} \leq \frac{\widetilde\Psi_{c,\text{min}}}{4} \left( 2 - \frac{4}{\widetilde\Psi_{c,\text{min}}} C_{D,2} \right)^{\frac{1}{q}}.
\label{eq:hyper_given_laplace}
\end{equation}

We set the Laplacian coefficient to an arbitrary stable value of $C_{D, 2}$ = 0.05, then seek the largest stable coefficient for the additional hyperdiffusion. The resulting maximum $C_{D, 2q}$ values at a C192 resolution in Table \ref{table:mixed_order_div_limits} show that, similarly to a single order of diffusion, the virtually-inviscid limits are slightly smaller than for default CAM transport, with both being close to the linear theory. The equiangular virtually-inviscid case that blew up with solely eighth-order damping (Table \ref{table:div_damp_limits}) is now stabilised with the addition of Laplacian diffusion. \par

\begin{table}[htpb]
\caption[]{JW2006 tests with combined Laplacian and hyperviscous divergence damping in CAM-FV3. The linear stability limits (the first row of the values) are calculated using (\ref{eq:hyper_given_laplace}) and are rounded down to three decimal places. A C192 resolution is used with two transport schemes (default CAM, virtually-inviscid). The Laplacian diffusion is set to $C_{D,2} = 0.05$ and the largest stable $C_{D,2q}$ is identified.}
\label{table:mixed_order_div_limits}
    \begin{center}
    \begin{tabular}{ |c||c|c|c||c|c|c| }
    \hline
         & \multicolumn{3}{|c||}{Equi-edge} & \multicolumn{3}{|c|}{Equiangular} \\
         \hline
         \hline
         Additional hyperviscosity order & 4th & 6th & 8th & 4th & 6th & 8th \\
        \hline
        \hline
        Linear stability limit & 0.185 & 0.170 & 0.163 & 0.147 & 0.137 & 0.132 \\
        \hline
        \hline
        Default CAM (monotonic) & 0.185 & 0.171 & 0.164 & 0.150 & 0.139 & 0.134 \\
        \hline
        Virtually-inviscid unlimited & 0.181 & 0.168 & 0.162 & 0.148 & 0.138 & 0.133 \\
        \hline
    \end{tabular}
    \end{center}
\end{table}

FV3 has an additional mixed-order divergence damping option, where the Laplacian diffusion coefficient varies in space whilst the hyperviscosity is fixed. This option is referred to as `Smagorinsky-type' diffusion in FV3 literature \parencite{harris2021scientific_FV3_doc}, but we instead call this \textit{flow-dependent divergence damping} to avoid confusion with the diffusion mechanism of \textcite{smagorinsky1963general}. Typical Smagorinsky damping sets the diffusion coefficient relative to the local deformation of $\sqrt{\epsilon^2 + \gamma^2}$, using quantities $\epsilon = u_x - v_y$ and $\gamma = u_y + v_x$ where subscripts denote partial differentiation, and is used in models such as MPAS \parencite{skamarock2012multiscale} and ICON \parencite{zangl2015icon}.  \par
FV3's flow-dependent Laplacian coefficient scales with the L2 norm of the horizontal divergence and relative vorticity, with $\nu_{D,2}$ (\ref{eq:nonD_damp_coeff}) redefined as
\begin{equation}
    \nu_{D,2} = C^* \Delta A_{\text{min}} \sqrt{D^2 + \zeta^2},
\label{eq:smag}
\end{equation}

\noindent with $C^*$ the flow-dependent diffusion coefficient. CAM-FV3 does not use flow-dependent diffusion by default, but when using this mechanism, \textcite{harris2021scientific_FV3_doc} suggests a value of 0.2 or 0.5. We set $C^* = 0.5$, then identify the largest stable $C_{D,2q}$ for hyperviscosity that can be used with the modified Laplacian damping. Whilst flow-dependent damping does reduce the maximum allowable hyperviscosity strength, as shown in Table \ref{table:div_damp_limits}, it is by a much smaller amount than with a spatially-constant Laplacian coefficient of $C_{D,2} = 0.05$ (Table \ref{table:mixed_order_div_limits}). For reference, at 850 hPa and day 15 in the JW2006 C96 test, $\max \{\sqrt{D^2+\zeta^2}\} \approx 2.5 \times 10^{-4} ~\text{s}^{-1}$. The flow-dependent diffusion will have more impact in flows with a larger factor of $\sqrt{D^2+\zeta^2}$, which then more tightly restricts the hyperviscosity coefficient for stability. \par

\subsection{Vorticity damping}
Vorticity damping is implemented in FV3 through diffusive fluxes \parencite{harris2021scientific_FV3_doc}. The diffused vorticity field, which is the action of the pseudo-Laplacian or higher-order operator on $\zeta$ on the primary D-grid, is
\begin{equation}
    \widetilde\nabla^{2q} \zeta = \frac{1}{\Delta A_d} \left[ \delta_x \lbrac \frac{\Delta y_d}{\Delta x_c} \sin(\alpha) \delta_x (\widetilde\nabla^{2(q-1)} \zeta) \rbrac + \delta_y \lbrac \frac{\Delta x_d}{\Delta y_c} \sin(\alpha) \delta_y (\widetilde\nabla^{2(q-1)} \zeta) \rbrac  \right].
    \label{eq:FV3_vort_higher_operator}
\end{equation}

In FV3, there are two terms in the horizontal momentum equations (\ref{eq:fv3_momentum_eqs}) that damp the vorticity. First are the $\mathcal{T}_X,\mathcal{T}_Y$ transport operators of \textcite{lin1996multidimensional} that contain implicit diffusion. Second are diffusive fluxes of $\mathcal{V}_x,\mathcal{V}_y$ that apply the optional explicit vorticity damping. This leads to a major difference between practical stability limits for divergence and vorticity damping in FV3: vorticity is implicitly diffused by the transport operators, whereas divergence is not directly damped during transport. Hence, whilst the maximum divergence damping diffusion coefficients in idealised testing are very close to linear theory (Table \ref{table:div_damp_limits}), the practical limits on $C_{\zeta,2q}$ are expected to be much lower, as the combination of both implicit and explicit vorticity damping must remain stable during each time step. \par
CAM-FV3 does not use vorticity damping by default, instead relying upon sufficient implicit diffusion from the transport scheme. Only fourth- and sixth-order vorticity damping are available, and this choice is determined from the order of the divergence damping --- fourth-order vorticity damping is used with Laplacian or fourth-order divergence damping, and sixth-order vorticity damping with sixth- or eighth-order divergence damping. Our tests will examine the fourth- and sixth-order mechanisms by setting the divergence damping to the same order, $q \in \{ 2,3 \}$. The default CAM divergence damping coefficient of $C_{D,2q}$ = 0.15 is retained for the equi-edge grid. As this value is unstable for sixth-order divergence damping on the equiangular grid (Fig. \ref{fig:equi_ang_scale_select}), a reduced $C_{D,2q}$ = 0.13 is used for this gnomonic mapping. \par
We again compare two unlimited and three monotonic transport schemes (see the appendix for more details). Typically, monotonic schemes are more diffusive due to the additional monotonicity constraints. The implicit diffusion strength depends on the form of the constraints; for example, the conditions for the Huynh monotonic scheme \parencite{huynh2007schemes} are more comprehensive than the Lin monotonic scheme \parencite{lin2004vertically}, which makes the Huynh scheme more computationally expensive, but likely less diffusive. \par
Table \ref{table:vort_limits_jw2006} verifies that the practical vorticity damping limits in the CAM-FV3 tests are lower than predicted by linear theory for purely explicit damping. In the JW2006 case, the practical limit ranges from 43\% of the theoretical limit for Lin monotonic transport on the C192 equi-edge grid and fourth-order vorticity damping, to 76\% with sixth-order vorticity damping on the C96 equiangular grid. Stronger vorticity damping can be used in the HS1994 tests (Table \ref{table:vort_limits_hs1994}), which indicates less implicit transport diffusion compared to the JW2006 tests. Intriguingly, a larger coefficient can be used with sixth-order vorticity damping compared to fourth-order; this contradicts the expectation from linear theory that higher-order explicit damping has a stricter stability limit (\ref{eq:damp_stab_limit}). For most cases, stronger vorticity damping can be applied on the equiangular grid than on the equi-edge grid.  \par
In the equi-edge JW2006 tests, the choice of horizontal transport scheme clearly affects the upper bounds on vorticity damping. In order of increasing maximum $C_{\zeta,2q}$, we have: Lin monotonic, intermediate unlimited, default CAM (monotonic), virtually-inviscid unlimited, and Huynh monotonic. When considering linear stability of the combination of implicit and explicit vorticity diffusion, we postulate that this gives a relative ranking of implicit transport diffusion in this specific test. A surprising result is that the Huynh monotonic scheme allows for stronger vorticity damping than the unlimited schemes in this case. Much smaller differences were found in the HS1994 tests, and none were observed in the equiangular grid JW2006 tests. \par

\begin{table}[htpb]
\caption[]{The maximum coefficient for vorticity damping, to three decimal places, that can be used in the JW2006 test case for fifteen days. The linear stability limits are computed from (\ref{eq:damp_stab_limit}) and rounded down to three decimal places. Fourth- and sixth-order vorticity damping are tested with five different transport schemes. These schemes are listed in order of increasing maximum stable vorticity damping coefficient on the equi-edge grid.}
\label{table:vort_limits_jw2006}\begin{center}
    \begin{tabular}{ |c||c|c||c|c||c|c||c|c| }
    \hline
         & \multicolumn{4}{|c||}{Equi-edge} & \multicolumn{4}{|c|}{Equiangular} \\
         \hline
         Grid resolution & \multicolumn{2}{|c||}{C96} & \multicolumn{2}{|c||}{C192} & \multicolumn{2}{|c||}{C96} & \multicolumn{2}{|c|}{C192} \\
         \hline
        Order of vorticity damping & 4th & 6th & 4th & 6th & 4th & 6th & 4th & 6th \\
        \hline
        Linear stability limit & 0.203 & 0.181 & 0.203 & 0.181  & 0.167 & 0.148 & 0.166 & 0.148 \\
        \hline
        \hline
        \hline
        Lin monotonic & 0.097 & 0.113 & 0.087 & 0.104 & 0.110 & 0.113 & 0.107 & 0.111 \\
        \hline
        Intermediate unlimited & 0.098 & 0.113 & 0.089 & 0.106 & 0.110 & 0.113 & 0.107 & 0.111 \\
        \hline
        Default CAM (monotonic) & 0.099 & 0.114 & 0.090 & 0.107 & 0.110 & 0.113 & 0.107 & 0.111 \\
        \hline
        Virtually-inviscid unlimited & 0.104 & 0.116 & 0.094 & 0.108 & 0.110 & 0.113 & 0.107 & 0.111 \\
        \hline
        Huynh monotonic & 0.105 & 0.119 & 0.098 & 0.113 & 0.110 & 0.113 & 0.107 & 0.111 \\
        \hline
    \end{tabular}
    \end{center}
\end{table}

\begin{table}[htpb]
\caption[]{The maximum coefficient for vorticity damping, to three decimal places, that can be used in the spun-up C96 HS1994 test for thirty days. The linear stability limits are computed from (\ref{eq:damp_stab_limit}) and rounded down to three decimal places.}
\label{table:vort_limits_hs1994}
\begin{center}
    \begin{tabular}{ |c||c|c||c|c|| }
    \hline
         & \multicolumn{2}{|c||}{Equi-edge} & \multicolumn{2}{|c|}{Equiangular} \\
         \hline
        Order of vorticity damping & 4th & 6th & 4th & 6th  \\
        \hline
        Linear stability limit & 0.203 & 0.181 & 0.167 & 0.148 \\
        \hline
        \hline
        \hline
        Lin monotonic & 0.110 & 0.122 & 0.114 & 0.116 \\
        \hline
        Intermediate unlimited & 0.111 & 0.123 & 0.114 & 0.116 \\
        \hline
        Default CAM (monotonic) & 0.110 & 0.121 & 0.114 & 0.116 \\
        \hline
        Virtually-inviscid unlimited & 0.112 & 0.123 & 0.115 & 0.117 \\
        \hline
        Huynh monotonic & 0.111 & 0.123 & 0.114 & 0.117 \\
        \hline
    \end{tabular}
    \end{center}
\end{table}

\clearpage

\section{Summary and discussion}
\label{section:discussion}
This paper investigated the stability of divergence and vorticity damping on the equidistant, equiangular, and equi-edge gnomonic cubed-sphere grids. A von Neumann analysis was performed on damping operators that either use a simplified pseudo-Laplacian, as in the FV3 dynamical core, or the full curvilinear Laplacian. The resulting linear stability bounds depend on the gnomonic mapping through the cubed-sphere cell areas, aspect ratios, and nonorthogonality factor of $\sin(\alpha)$. These quantities are encapsulated in the grid stability function, $\widetilde \Psi$ for the pseudo-Laplacian and $\Psi$ for the full Laplacian, which is evaluated on the primary or offset grid depending on the locations of horizontal divergence and relative vorticity. The minimum value of $\widetilde \Psi$ or $\Psi$ dictates the theoretical linear stability limit. \par
When using the pseudo-Laplacian operator, the equiangular grid has the tightest restriction on the damping coefficient for linear stability, whilst the equidistant and equi-edge grids require a smaller coefficient when using the full Laplacian. For the equiangular grid, the location of diffusive instability is at the middle of panel edges with the pseudo-Laplacian, but moves to the panel corners when using the full Laplacian. The equidistant and equi-edge grids are most unstable at the panel corners for both damping operators. Comparisons of $2 \Delta x$ wave amplification factors across a cubed-sphere panel showed that a nondimensional diffusion coefficient using the minimum cell area can lead to large variations in damping strength. This damping inhomogeneity is exacerbated with higher-order hyperdiffusion and can lead to the grid-scale waves being minimally diffused in the largest cells. Of the three gnomonic mappings, the equiangular grid leads to the most uniform damping of the $2 \Delta x$ wave with both the pseudo-Laplacian and full Laplacian operators. \par
The D-grid CAM-FV3 model was used to practically examine the analytical stability limits from the pseudo-Laplacian operator. The CAM-FV3 default divergence damping coefficient of $C_{D,2q}=0.15$, which is tuned for the equi-edge grid, is unstable for sixth- and eighth-order divergence damping on the equiangular grid. The default coefficient also leads to negative amplification factors for the small-scale waves on both grids, so an alternative oscillation-free coefficient was proposed, which completely damps the $2\Delta x$ wave and ensures positive amplification factors for all other waves. It is worth investigating whether the stricter oscillation-free coefficient is advantageous to use in dynamical cores, or if potential oscillations introduced by diffusion have no practical bearing on model performance. \par
Divergence damping tests with CAM-FV3 in baroclinic wave and Held-Suarez simulations showed good agreement with linear theory. Some tests with eighth-order divergence damping were unconditionally unstable, and either adding Laplacian damping or switching to sixth-order damping helped stabilise these simulations. Practical limits to vorticity damping were well below the linear stability limits due to implicit vorticity diffusion in the FV3 transport operators. Differences in the upper bound of $C_{\zeta,2q}$ between horizontal transport schemes were most clearly observed for the baroclinic wave test on the equi-edge grid. Further numerical studies would be useful to ascertain the efficacy of using vorticity damping stability as a measure of implicit diffusion in the FV3 transport scheme. It is also worth making comparisons to other approaches of estimating the implicit diffusion. \par
As an alternative to stability bounds derived from the von Neumann method, a future study could investigate eigenvalues of the diffusion operator, similarly to the work of \textcite{cheong2015eigensolutions} for a spectral element discretisation. It would also be useful to examine the impact on overall model behaviour of damping with a pseudo-Laplacian instead of the full curvilinear Laplacian. Additionally, studies into the lower bound on the divergence and vorticity damping coefficients would be valuable to complement this study of the upper bound.

\section*{Acknowledgments}
This work was funded by the National Oceanic and Atmospheric Administration (NOAA), grant NA22OAR4320150-T3-01S093.

\section*{Code Availability}
A clone of the publicly available Community Atmosphere Model code repository (\url{https://github.com/ESCOMP/CAM}), set to tag \texttt{6\_4\_050}, was used for the CAM-FV3 tests in this work. A collection of analysis scripts and commands for the CAM-FV3 simulations can be found in the GitHub repository \\ \url{https://github.com/ta440/div_vort_damp_cubed_sphere}.

\printbibliography

\end{document}